\documentclass[11pt]{article}

\usepackage[left=1in,top=1in,right=1in,bottom=1in,nohead]{geometry}

\usepackage{latexsym,amsmath,amssymb, amsfonts,enumerate,
  epsfig, graphics,times,color, amsthm,url}

\newtheorem{theorem}{Theorem}
\theoremstyle{definition}
\newtheorem{lemma}[theorem]{Lemma}
\newtheorem{remark}{Remark}
\newtheorem{algorithm}{Algorithm}

   \newtheoremstyle{example}{\topsep}{\topsep}%
     {\it}
     {}
     {\bfseries}
     {}
     {\newline}
     {\thmname{#1}\thmnumber{ #2}\thmnote{ #3}}

   \theoremstyle{example}
   \newtheorem{example}{Example}

\newcommand{\eqdef}{\overset{\mbox{\tiny def}}{=}}
\newcommand{\R}{\mathbb{R}}
\newcommand{\Z}{\mathbb{Z}}
\newcommand{\E}{\mathbb{E}}
\newcommand{\eqdist}{\overset{\mathcal{D}}{=}}

\newcommand{\X}{X^N}
\newcommand{\Zm}{Z^N}
\newcommand{\hl}{h_{\ell}}

\newcommand{\Ql}{\widehat{Q}}
\newcommand{\N}{\overline N}

\title{Multi-level Monte Carlo for continuous time Markov chains,\\
        with applications in biochemical kinetics}

\author{David F. Anderson$^{1}$ and Desmond J. Higham$^{2}$}

\begin{document}

\maketitle

\begin{abstract}
We show how to extend a recently proposed multi-level Monte Carlo approach 
to the continuous time Markov chain setting, thereby greatly lowering the computational complexity needed to compute expected values of functions of the state of the system to a specified accuracy.  The extension is non-trivial, exploiting a coupling of the requisite processes that is  
 easy to simulate while providing a  small variance for the estimator.  Further, and in a stark departure from  other  implementations of multi-level Monte Carlo, we show how to 
produce an unbiased estimator that is significantly less computationally expensive than the usual unbiased estimator arising from exact algorithms in conjunction with  crude Monte Carlo.   We thereby dramatically improve,
in a quantifiable manner,
the basic computational complexity of current approaches that
have many names and variants across the scientific literature,
including the Bortz-Kalos-Lebowitz algorithm,
discrete event simulation,
dynamic Monte Carlo,
kinetic Monte Carlo,
the n-fold way,
the next reaction method,
the residence-time algorithm,
the stochastic simulation algorithm,
Gillespie's algorithm,
and
tau-leaping.
The new algorithm applies generically, but we also give an example where the
coupling idea alone, even without a multi-level discretization, can be used to improve efficiency
by exploiting system structure.
Stochastically modeled chemical reaction networks provide a very
important application for this work. Hence, we use this context for
our notation, terminology, natural scalings, and computational examples.


\end{abstract}

\footnotetext[1]{Department of Mathematics, University of
  Wisconsin, Madison, Wi. 53706, anderson@math.wisc.edu, grant support from NSF-DMS-1009275.}

\footnotetext[2]{Department of Mathematics and Statistics, University of
  Strathclyde, Glasgow, G1 1XH, d.j.higham@strath.ac.uk, supported by 
 a Fellowship from the Leverhulme Trust.}



\section{Introduction}
\label{sec:intro}

This paper concerns the efficient computation of expectations for continuous time Markov chains.  Specifically, we extend  the multi-level Monte Carlo approach
of Giles \cite{Giles2008}, with related earlier work 
by Heinrich 
\cite{Heinrich2001}, to this setting.  We study the wide class of systems that can be written using the  random time change representation of Kurtz \cite[Chapter 6]{Kurtz86} \cite{KurtzPop81} in the 
form 
\begin{align}
  X(t) &= X(0) + \sum_{k=1}^R Y_k\left( \int_0^t \lambda_k(X(s))ds
  \right)\zeta_k
   \label{eq:RTC_1st},
\end{align}
where the $Y_k$ are independent unit-rate Poisson processes,  $\zeta_k \in \R^d$, and the functions $\lambda_k$ are the associated intensity, or propensity, functions.  
While such models are used in nearly all branches of the sciences, especially in the studies of  queues and  populations, their use has recently exploded in the biosciences, and we use this application area for the setting of our work.  We will formally introduce these models  in Section \ref{sec:model}, however we begin by demonstrating how two different models, one from chemistry and one from queuing, can be represented via \eqref{eq:RTC_1st}.  

First, consider a linear reversible chemical network 
\begin{equation*}
   S_1 \overset{\kappa_1}{\underset{\kappa_2}{\rightleftarrows}} S_2,
\end{equation*} 
in which molecules of type $S_1$ convert to molecules of type $S_2$ at rate $\kappa_1 X_1$, where $X_1$ is the number of $S_1$ molecules, and molecules of type $S_2$ convert to $S_1$ at rate $\kappa_2 X_2$.  Here we are assuming the system satisfies mass action kinetics, see Section \ref{sec:model}.  The usual stochastic model, written in the framework of \eqref{eq:RTC_1st}, is then
\begin{align*}
	X(t) &=X(0) + Y_1\left(\int_0^t \kappa_1 X_1(s) ds\right)\left(\begin{array}{c}
	-1\\
	1
	\end{array}\right) +  Y_2\left(\int_0^t \kappa_2 X_2(s) ds\right)\left(\begin{array}{c}
	1\\
	-1
	\end{array}\right). 
\end{align*}

Next, consider an $M/M/k$ queue in which arrivals are happening at a constant rate $\lambda >0$, and there are $k$ servers, with each serving at a rate $\mu>0$.  Letting $X(t)$ denote the number of customers in the queue at time $t$,
\begin{equation*}
	X(t) = X(0) + Y_1\left(\lambda t\right) - Y_2\left( \mu \int_0^t (X(s) \wedge k ) \ ds) \right),
\end{equation*}
where we define $a\wedge b \eqdef \min\{a,b\}$.


There  are multiple algorithms available to compute exact sample paths of  continuous time Markov chains, and, though they are all only slight variants of each other, they go by different names depending upon the branch of science within which they are being applied.  These  include the Bortz-Kalos-Lebowitz algorithm,
discrete event simulation,
dynamic Monte Carlo,
kinetic Monte Carlo,
the n-fold way,
the residence-time algorithm,
 the stochastic simulation algorithm, the next reaction method, and Gillespie's algorithm, where the final two are the most commonly referred to algorithms in the biosciences.
As the computational cost of exact algorithms scales linearly with the number of jump events (i.e. reactions), such methods become computationally intense for even moderately sized systems.  This 
issue looms large when many sample paths are needed in a Monte Carlo setting.     
To address this, approximate algorithms, and notably the class of algorithms termed ``tau-leaping'' methods introduced by Gillespie \cite{Gill2001} in the chemical kinetic setting, have been developed with the explicit aim of greatly lowering the computational complexity of each path simulation while controlling the
bias \cite{Anderson2007b,AndersonGangulyKurtz,AndMaso2011,Li2011,Li2007,Rathinam2005}. 

A common task in the study of stochastic models, 
and the main focus of this paper,  is to approximate $\E f(X(T))$, where $f$ is a scalar-valued 
function of the state of the system which gives a measurement of interest.    For example, the function $f$ could be:
\begin{enumerate}
\item $f(X(T)) = X_i(T)$, yielding estimates for mean values, or 
\item $f(X(T)) = X_i(T)X_j(t)$,  which can be used with estimates for  the mean values to provide estimates of variances (when $i = j$) and covariances (when $i \ne j$), or
\item $f(X(T)) = 1_{\{X(T) \in B\}}$,  the indicator function giving 1 if the state is in some specified set.  Such functions could also be used to construct histograms, for example, since $\E f(X(T)) = P\{X(T) \in B\}$.
\end{enumerate}

Suppose we use an exact simulation algorithm to approximate $\E f(X(T))$  to 
$O(\epsilon)$ accuracy 
in the sense of confidence intervals.  To do so, we  need to generate $n=O(\epsilon^{-2})$ paths so that the standard deviation of the usual Monte Carlo estimator,
\begin{equation*}
	\mu_n = \frac{1}{n} \sum_{j = 1}^n f(X_{[j]}(T)),
\end{equation*} 
where $X_{[j]}$ are independent realizations generated via an exact algorithm, is $O(\epsilon)$.  If we let $\overline N>0$ be the order of magnitude of the number of computations needed to produce a single sample path using an exact algorithm, then the total computational complexity becomes  $O(\overline N \epsilon^{-2})$.  (Here, and throughout, we work in terms of expected 
computational complexity.)

When $\N \gg 1$, which is the norm as opposed to the exception in many settings, it may be desirable to make use of an approximate algorithm.  Suppose $\E f(X(T)) - \E f(Z_h(T)) = O(h)$, where $Z_h$ is an approximate path generated from a time discretization with a magnitude of $h$ 
(i.e.\ we have a weakly order one method).  We first make the trivial observation that  the estimator 
\begin{equation}
	\mu_n = \frac{1}{n} \sum_{j = 1}^n f(Z_{h,[j]}(T)),
	\label{eq:crude_approx}
\end{equation}
where $Z_{h,[j]}$ are independent paths generated via the approximate algorithm with a step size of $h$, is an unbiased estimator of $\E f(Z_{h}(T))$, and not $\E f(X(T))$.  However, noting that 
\begin{equation}
	\E f(X(T)) -\mu_n= \big[\E f(X(T)) - \E f(Z_{h}(T))\big] +\big[ \E f(Z_{h}(T)) - \mu_n\big],
	\label{eq:bias_sd}
\end{equation}
we see that choosing $h = O(\epsilon)$, so that the first term on the right is $O(\epsilon)$, and  $n = O(\epsilon^{-2})$, so that the standard deviation is $O(\epsilon)$,  delivers the desired accuracy.
With a fixed cost per time step, the computational complexity of generating a single such path is $O(\epsilon^{-1})$ and we find that the total computational complexity is $O(\epsilon^{-3})$.  Second order methods lower the computational complexity to $O(\epsilon^{-2.5})$, as $h$ may be chosen to be $O(\epsilon^{1/2})$. 

The discussion above suggests that the choice between exact or approximate path computation
should depend upon  
whether $\epsilon^{-1}$ or $\overline N$ is the larger value, with an exact algorithm being 
beneficial when $\overline N < \epsilon^{-1}$.   It is again worth noting, however, that the estimators built from approximate methods are biased, and while analytic bounds can be provided for that bias \cite{AndersonGangulyKurtz,AndMaso2011,Li2007} these are
typically neither sharp nor computable, 
and hence of limited practical value.
The exact algorithm, on the other hand,  trivially produces an {\em unbiased} estimator, so 
it may be argued that $\epsilon^{-1}\ll \overline N$ is necessary before it is worthwhile to 
switch to an approximate method.

In the diffusive setting
the multi-level Monte Carlo approach
 has the remarkable property of 
lowering the standard 
 $O(\epsilon^{-3})$
cost of computing an 
 $O(\epsilon)$ accurate Monte Carlo estimate of 
$\E f(X(T))$ 
down to  $O(\epsilon^{-2}\log(\epsilon)^2$) \cite{Giles2008}. Here, 
we are assuming 
that a weak order one and strong order $1/2$ discretization method, such as Euler--Maruyama, is used. Further refinements have appeared in
\cite{Giles2009, GHM2009, Higham2011, Kloeden2011,Kloeden2011b},
 and the same ideas have been 
applied to partial differential equations \cite{Barth2011,Cliffe2011}. 
  A key motivation for multi-level Monte Carlo is that optimizing the overall expected value computation is a different, and typically more relevant, goal than optimizing along each path. 
Computing an expectation using only an exact algorithm (or an algorithm with a very fine time-step) can require a large number of paths and an {\em extremely} large number of random variables and state updates. In general, the total number of paths cannot be reduced.
 The
 computational benefits of multi-level Monte Carlo arise because the 
 number of random variables and state updates needed to approximate
 the expectation can be drastically reduced by averaging over a very carefully chosen combination of coordinated realizations, many of which are much cheaper to compute than an exact realization.

In this paper we extend the multi-level approach to the continuous time Markov chain setting, and especially the stochastic chemical kinetic setting.  The extension involves a non-trivial
  coupling of the requisite processes that is easy to simulate while providing a very small variance for the estimator.  In fact, showing the
practical importance of the coupling (found in this paper in both equations \eqref{eq:unscaled_coupled_tau} and \eqref{eq:unscaled_coupled_exact}), which was 
first used in \cite{Kurtz82} and later in  \cite{AndersonGangulyKurtz} as an analytical tool and  subsequently in \cite{AndersonSensitivities} towards the problem of computing parameter sensitivities,  could be viewed as the most important contribution of this paper.   Further,  and in a stark departure from other  implementations of multi-level Monte Carlo, we provide a second multi-level Monte Carlo algorithm which exploits the representation \eqref{eq:RTC_1st} to produce an {\em unbiased} estimator giving the desired accuracy with significantly less computational complexity than an exact algorithm alone.  The authors believe that this unbiased multi-level Monte Carlo will become a standard, generic algorithm for approximating expected values of continuous time Markov chains, and especially stochastically modeled chemical reaction networks.
 

We emphasize that the gains in computational efficiency
reported in this work apply to generic models, and do not rely on any specific structural properties.
However, the ideas have the potential to be fine-tuned further in appropriate cases; for example
by exploiting known analytical results or multi-scale partitions.   We provide such an example in Section \ref{sec:examples}.  We also emphasize that our complexity analysis does not involve asymptotic limits.  In particular, we do not consider infinitely large system size, where stochastic effects vanish, or infinitesimally small discretization time-step, where the benefits of an approximate method evaporate.

 The outline of the remainder of the paper is as follows.  In Section \ref{sec:model}, we 
consider stochastically modeled chemical reaction networks, which is our main application area, discussing how such models can be represented via \eqref{eq:RTC_1st}.  In Section \ref{sec:scalings}, we introduce an equivalent  model to \eqref{eq:RTC_1st} that incorporates the natural temporal and other quantitative scales. 
Consideration of such a scaled model is critical for realistic quantitative comparisons
 of accuracy versus cost for computational methods, though it plays no role in the actual simulations.  In Section \ref{sec:tau-leaping}, we briefly review Euler's method, often called tau-leaping in the chemical kinetic setting.  In Section \ref{sec:old_methods}, we review the original multi-level Monte Carlo method.  In Section \ref{sec:new_methods}, we extend multi-level Monte Carlo to the continuous time Markov chain setting in two different ways.  In the first, exact algorithms are not used and we are led to an efficient method with an unquantified bias.  In the second, exact algorithms play a key role and allow us to develop  \textit{unbiased} estimators. 
In both cases,
we quantify precisely the generic computational efficiencies obtained, relative to standard 
Monte Carlo. In Section \ref{sec:analysis}, we provide the delayed proofs of the main analytical results of Section \ref{sec:new_methods}.   In Section \ref{sec:implementation}, we briefly discuss some implementation issues.  In Section \ref{sec:examples}, we provide  computational examples demonstrating our main results.  Finally, in Section \ref{sec:conclusions} we provide some brief conclusions.

\section{The basic stochastic  model for chemical reaction networks}
\label{sec:model}

In this section we discuss how the basic stochastic model for chemical reaction networks can be represented via \eqref{eq:RTC_1st} for suitable choices of $\lambda_k$ and $\zeta_k$.  
A chemical reaction network consists of the interaction of multiple {\em species}, $\{S_1,\dots,S_d\}$,  through different possible {\em reactions}.  If we denote by $\zeta_k \in \R^d$ the change to the state of the system after each occurrence of the $k$th reaction, then we have
\begin{equation*}
	X(t) = X(0) + \sum_k R_k(t) \zeta_k,
\end{equation*}
where $X_i(t)$ gives the number of molecules of $S_i$ at time $t$, and $R_k(t)$ is the number of times the $k$th reaction has taken place up until time $t$.  
To model $R_k$, each reaction channel is assumed 
to have an associated intensity, or propensity, function, $\lambda_k:\R^d \to \R_{\ge 0}$, and for the standard Markov chain model, the number of times that the $k$th reaction occurs by time $t$ can then be
represented by the counting process 
\begin{equation*}
 R_k(t) = Y_k\bigg(\int_0^t \lambda_k(X(s))ds\bigg), 
 \end{equation*}
  where the $Y_k$ are independent
unit-rate Poisson processes; see, for example, \cite{KurtzPop81}, \cite[Chapter~6]{Kurtz86}, or the recent survey \cite{AndKurtz2011}.  The state of the system then satisfies 
\eqref{eq:RTC_1st}.
This formulation is termed a ``random time change representation''
and is equivalent to the ``chemical master equation representation''
found in much of the biology and chemistry literature. 


A common choice of intensity function for chemical reaction systems,
and the one we adopt throughout, is that of {\em mass action kinetics}.  Under mass
action kinetics, the intensity function for the $k$th 
reaction 
is
\begin{equation}
\lambda_k(x) = \kappa_k \prod_{i = 1}^d \frac{x_i!}{(x_i - \nu_{ki})!}1_{\{x_i \ge \nu_{ki}\}},
  \label{eq:stoch_MA}
\end{equation}
where $\nu_{ki}$ denotes the number of molecules of $S_i$ required for one instance of the reaction.  
Note that  $\lambda_k(x) = 0$ whenever $x_i \le 0$ and $\nu_{ki} \ne 0$.  We note that none of the core ideas of this paper depend upon the fact that $\lambda_k$ are mass-action kinetics and the assumption is made for analytical convenience
and historical
 consistency.

This  model is a
continuous time Markov chain in $\Z^d$ with generator
\begin{equation*}
  ({\mathcal A} f)(x) = \sum_{k} \lambda_k(x)(f(x + \zeta_k) - f(x)),
\end{equation*}
where $f : \Z^d \to \R$.  Kolmogorov's forward equation, termed the {\em chemical master equation} in much of the biology literature,
for this model is
\begin{equation*}
  \frac{d}{dt} P(x, t | \pi) = \sum_k \lambda_k(x-\zeta_k) 1_{\{x - \zeta_k \in \Z^d_{\ge 0}\}} P(x-
  \zeta_k,t|\pi) - \sum_k \lambda_k(x) P(x,t|\pi), 
\end{equation*}
where for $x \in \Z^d_{\ge 0}$, $P(x,t|\pi)$ represents the probability
that $X(t) = x$, conditioned upon the initial distribution $\pi$.

\begin{example}
  To solidify notation, we consider the network
  \begin{equation*}
  	S_1 \overset{\kappa_1}{\underset{\kappa_2}{\rightleftarrows}} S_2,\qquad 2S_2 \overset{\kappa_3}{\rightarrow} S_3,
  \end{equation*}
  where we have placed the rate constants $\kappa_k$ above or below their respective reactions.  For this example, equation \eqref{eq:RTC_1st} is
  \begin{align*}
  	X(t) = X(0) &+ Y_1\left( \int_0^t \kappa_1 X_1(s)ds\right)\left[\begin{array}{c}
	-1\\
	1\\
	0
	\end{array} \right] + Y_2\left( \int_0^t \kappa_2 X_2(s)ds\right)\left[\begin{array}{c}
	1\\
	-1\\
	0
	\end{array} \right] \\
	&+ Y_3\left( \int_0^t \kappa_3 X_2(s)(X_2(s)-1)ds\right)\left[\begin{array}{c}
	0\\
	-2\\
	1
	\end{array} \right].
  \end{align*}
 Using $\zeta_1 = [-1,1,0]^T$, $\zeta_2 = [1,-1,0]^T$, and $\zeta_3 = [0,-2,1]^T$, the generator ${\mathcal A}$ satisfies
  \begin{equation*}
  	({\mathcal A} f)(x) = \kappa_1 x_1(f(x +\zeta_1 ) - f(x)) + \kappa_2 x_2(f(x + \zeta_2) - f(x))+ \kappa_3 x_2(x_2 - 1)(f(x + \zeta_3) - f(x)). 
  \end{equation*}
\end{example}

\section{Scaled models}
\label{sec:scalings}

To quantify 
the relative computational complexity of different methods, it is important that the natural scalings of a model be taken into account.  However, we stress that such a change to the representation of the model does not change the simulation---we simulate the unscaled model but analyze
the methods on an appropriately scaled version.  

Letting $N$ be some natural parameter of the system, which is usually taken to be the abundance of the most abundant component, we scale the model by setting $X_i^N = N^{-\alpha_i}X_i$, where $\alpha_i \ge 0$ is chosen so that $X^N_i = O(1)$.    The general form of such a scaled model is 
\begin{equation}
	X^N(t) = X^N(0) + \sum_k Y_k\left( N^{\gamma} \int_0^t N^{c_k} \lambda_k(X^N(s)) ds\right)\zeta_k^N,
	\label{eq:main_multi}
\end{equation}
where $\gamma$ and $c_k$ are scalars, $|\zeta_k^N| = O(N^{-c_k})$, and both $X^N$ and $\lambda_k(X^N)$ are of order one.  Note that we are explicitly allowing for $|\zeta_k^N|$ to be smaller than $N^{-c_k}$, a point made explicit in and around equation \eqref{eq:rho_k}.   We note that we should write $\lambda_k^N$, as the resulting intensity function may depend upon $N$, though we drop the superscript $N$ for notational convenience.  It is now natural to take
\begin{equation*}
	\overline N = N^{\gamma}\sum_k N^{c_k}
\end{equation*}
as the order of magnitude for the number of computations required to generate a single path using an exact algorithm.  We will demonstrate how to arrive at such a scaled model for  chemical systems below, however we first discuss the parameter $\gamma$.

The parameter $\gamma$ should be interpreted as being related to the natural time-scale of the model.  That is, if $\gamma>0$ then the shortest timescale in the problem is much smaller than 1, while if $\gamma < 0$ it is much larger.  The analysis in this paper is  most applicable in the case that $\gamma \le 0$, 
for otherwise the error bounds grow quite rapidly.  
However, and as will be demonstrated in the examples section,  
the methods developed can still behave very well even 
when $\gamma > 0$, pointing out that the present analysis 
does not fully capture the behavior of the methods.

We will show how to derive a model of the form \eqref{eq:main_multi} in the case of chemical reaction networks with mass action kinetics.  Let $N \gg 1$, where  $N$ is  the abundance of the most abundant species, or some other large parameter.  Suppose we have a model of the form
\begin{equation*}
  X(t)=X(0)+\sum_kY_k\left(\int_0^t\lambda_k'(X(s))ds\right)\zeta_k,
\end{equation*} 
where the $\lambda_k'$ are of the form 
\begin{equation*}
  \lambda_k'(x)=\kappa'_k\prod_i \frac{x_i!}{(x_i - \nu_{ki})!}.
\end{equation*}
  For each species, define the {\em normalized
  abundance\/}  by
  $\X_i(t) \eqdef N^{-\alpha_i}X_i(t),$ 
where $\alpha_i\geq 0$ should be selected so that $\X_i=O(1)$. 
Here $X_i^N$  may be the species number ($\alpha_i=0$), the
species concentration, or something else.  Since the rate constants may also vary over several orders of
magnitude, we write $\kappa_k'=\kappa_kN^{\beta_k}$ where the
$\beta_k$ are selected so that $ \kappa_k=O(1)$.  Under  the mass-action kinetics assumption, we have that $\lambda_k'(X(s)) = N^{\beta_k + \nu_k\cdot \alpha} \lambda_k(X^N(s)),$ where $\lambda_k$ is deterministic mass-action kinetics with parameter $\kappa_k$ \cite{KangKurtz2011}, and we recall that $\nu_k$ is the source vector of the $k$th reaction.   
Our model has therefore become
\begin{equation*}  
  \X(t)= \X(0)+\sum_k Y_k\left(\int_0^tN^{\beta_k+\nu_
    k\cdot\alpha} \lambda_k(\X(s))ds\right)\zeta_k^N, 
\end{equation*}
where $\zeta_{ki}^N \eqdef \zeta_{ki}/N^{\alpha_i}$ (so $\zeta_k^N$ is the scaled reaction vector).  
Define $\gamma\in \R$ via
\begin{align*}
\begin{split}
\gamma \eqdef \max_{\{i,k \ : \ \zeta_{ki}^N \ne 0\} } \{ \beta_k + \nu_k \cdot \alpha - \alpha_i\}.
\end{split}
\end{align*}
 Then, for each $k$ define
\begin{equation}\label{eq:gamma}
c_k \eqdef \beta_k + \nu_k \cdot \alpha - \gamma.
\end{equation}
With these definitions, our chemical model becomes 
\eqref{eq:main_multi}.

\vspace{.2in}

Returning to the general setting of \eqref{eq:main_multi}, 
for each  $k$ we define
\begin{equation}\label{eq:rho_k}
	\rho_{k} \eqdef  \min\{ \alpha_i\ : \ \zeta^N_{ki} \ne 0 \},
\end{equation}
so that $|\zeta_k^N| \approx N^{-\rho_k}$, and define $\rho \eqdef \min\{\rho_k\}$.
We have that $\rho \ge 0$, and by the choice of $\gamma$ we have $c_k - \rho_k \le 0$ for all $k$.  Further, we point out that $\gamma$ is chosen so that $c_k = 0$
for at least one $k$.  Also, if $\|\nabla f \|_\infty$  is bounded, then 
\begin{equation*}
    N^{c_k} (f(x + \zeta_k^N) - f(x)) =  O(N^{c_k - \rho_k}),
\end{equation*}
with $c_k - \rho_k = 0$ for at least one $k$.  Finally, it is worth explicitly noting that the {\em classical scaling} holds if and only if $c_k \equiv \rho_k \equiv 1$ and $\gamma = 0$ \cite{AndersonGangulyKurtz,Kurtz72}.

%
%

\begin{remark}
We emphasise   that the models  \eqref{eq:RTC_1st} and \eqref{eq:main_multi} are equivalent in that $X^N$ is the scaled version of $X$.  The scaling is essentially an analytical tool as now both $X^N$ and $\lambda_k(X^N(\cdot))$ are $O(1)$, and in Section \ref{sec:analysis} it will be shown how the representation \eqref{eq:main_multi} is useful in the quantification of the behavior of different computational methods.  However, we stress that the scaling itself plays {\em no role} in the actual simulation of the processes, with the small exception that it can inform the decision for the size of the  time step of an approximate method.
\end{remark}

\begin{example}\label{ex:gamma}
To  solidify notation, consider the reversible isometry 
\begin{equation*}
	S_1 \overset{100}{\underset{100}{\rightleftarrows}} S_2
\end{equation*}
with $X_1(0) = X_{2}(0) = \text{10,000}$.  In this case, it is natural to take $N = $ 10,000 and $\alpha_1 = \alpha_2 = 1$.  As the rate constants are $100 = \sqrt{\text{10,000}}$, we take $\beta_1 = \beta_2 = 1/2$ and find that $\gamma = 1/2$ and $\rho_1 = \rho_2 = 1$.  The normalized process $X^N_1$ satisfies
\begin{equation*}
    X_1^N(t) = X_1^N(0) - Y_1\bigg(N^{1/2}N \int_0^t X_1^N(s)ds\bigg)\frac{1}{N} + Y_2\bigg(N^{1/2}N \int_0^t (2 - X_1^N(s))ds\bigg)\frac{1}{N}, 
\end{equation*}
where we have used that $X^N_1 + X^N_2 \equiv 2$.
\end{example}

\begin{example}
	We provide a deterministic example to further explain the use of the scalings.  Consider the ordinary differential equation
	\begin{equation*}
		\dot x(t) = \lambda N - \mu x(t),
	\end{equation*}	
	where $\lambda,\mu = O(1)$,  $N \gg 1$, and $x_0 = O(N)$.  Of course, the solution to this system is
	\begin{equation*}
		x(t) = \frac{\lambda N}{\mu}  - \left(\frac{\lambda N}{\mu} - x_0\right) e^{-\mu t}. 
	\end{equation*}	
	However, defining $x^N = N^{-1}x$, we see that $x^N$ satisfies
	\begin{equation*}
		\dot x^N(t) = \lambda - \mu x^N(t),
	\end{equation*}
	with $x^N_0 = O(1)$.  Solving yields
	\begin{equation*}
		x^N(t) = \frac{\lambda}{\mu} - \left(\frac{\lambda}{\mu} - x^N_0\right) e^{-\mu t}.
	\end{equation*}
	Note, then,  that solving for either $x$ or $x^N$ automatically yields the other after scaling.  Also note the important property that in the ODE governing $x$, the driving force, $\lambda N$, was an extremely large value.  However, the forcing function of $x^N$, which is simply $\lambda$, was $O(1)$.  
	\label{ex:ODE}
\end{example}

Example \ref{ex:ODE} points out an important feature:  the functions $\lambda_k$ of \eqref{eq:main_multi}, together with their derivatives, are much better behaved, in terms of their magnitude, than the intensity functions of the original model \eqref{eq:RTC_1st}.  
Therefore, after possibly redefining the kinetics by multiplication with a cutoff function, see, for example, \cite{AndersonGangulyKurtz, AndMaso2011}, it is reasonable to assume that each $\lambda_k$ is, in fact, a globally Lipschitz function of $X^N$.  We formalize this assumption here.

\vspace{.2in}

\noindent \textit{Running assumption:}  Throughout, we assume that the functions $\lambda_k$ of \eqref{eq:main_multi} are globally Lipschitz.\vspace{.2in}

\section{A review of  Euler's method in the current setting}
\label{sec:tau-leaping}

We briefly review Euler's method, termed tau-leaping in the chemical kinetic literature \cite{Gill2001}, as applied to the models \eqref{eq:RTC_1st}, and equivalently \eqref{eq:main_multi}.  
The basic idea of tau-leaping is
to hold the intensity functions fixed over a time interval $[t_n,
t_n+h]$ at the values $\lambda_k(X(t_n))$, where $X(t_n)$ is the
current state of the system, and, under this assumption, compute the
number of times each reaction takes place over this period.  
 As the
waiting times for the reactions are exponentially distributed this
leads to the following algorithm, which simulates up to a time of $T>0$.  
Below and in the sequel, for $x \ge 0$ we will write 
Poisson$(x)$ to denote a sample from the Poisson 
distribution with parameter $x$, with all
such samples 
being independent of each other and of all 
other sources of randomness used.

\begin{algorithm}[Euler tau-leaping]   \label{alg:Euler}
 Fix $h>0$.  Set $Z_h (0) = x_0$, $t_0 = 0$,
  $n=0$ and repeat the following until $t_{n} = T$:
  \begin{enumerate}[$(i)$]
   \item Set $t_{n+1} = t_{n} + h$.  If $t_{n+1} \ge T$, set $t_{n+1} = T$ and $h = T - t_n$.
  \item For each $k$, let $\Lambda_k = \text{Poisson}(\lambda_k(Z_h(t_n))h)$.
  \item Set $Z_h(t_{n+1}) = Z_h(t_{n}) + \sum_k \Lambda_k \zeta_k$.
   \item Set $n \leftarrow n+1$.
  \end{enumerate}
\end{algorithm}

 Several improvements and modifications have been made to the basic
algorithm described above over the years.
Some concern adaptive step-size selection 
along a path 
\cite{Cao2006, Gill2003}. Others focus on 
ensuring non-negative population values
\cite{Anderson2007b, Cao2005, Chatterjee2005, Tian2004}.
The latter issue 
is easily dealt with in our context; 
for example, 
it is sufficient to 
return a value to zero if it ever goes negative in the course 
of a simulation.  
This is discussed further in subsection~\ref{subsec:obs}.

Analogously 
to \eqref{eq:RTC_1st}, a path-wise representation of Euler tau-leaping defined for all $t\ge 0$
can be given through a random time change of Poisson processes:
\begin{equation}
  Z_h(t) = Z_h(0) + \sum_k Y_k \left( \int_0^t \lambda_k(Z_h \circ \eta(s))
    ds  \right)\zeta_k, 
  \label{eq:RTC_tau}
\end{equation}
where the $Y_k$ are as before, and $\displaystyle \eta(s) \eqdef \left \lfloor \frac{s}{h} \right \rfloor h$. Thus, $Z_h\circ \eta(s) = Z_{h}(t_n)$ if $t_n\le s < t_{n+1}$.  Noting that 
\begin{equation*}
  \int_0^{t_{n+1}} \lambda_k(Z_h \circ \eta(s)) ds = \sum_{i=0}^n \lambda_k(Z_h(t_i))(t_{i+1} - t_i)
  \end{equation*}
  explains why this method is called Euler tau-leaping.   
Following \eqref{eq:main_multi}, for each $i\in \{1,\dots,d\}$ we let $Z_{h,i}^N \eqdef N^{-\alpha_i} Z_{h,i}$,  so the scaled version of \eqref{eq:RTC_tau} is
\begin{equation}
	Z_h^N(t) = Z_h^N(0) + \sum_k Y_k \left( N^{\gamma} \int_0^t N^{c_k} \lambda_k(Z^N_h \circ \eta(s))
    ds  \right)\zeta_k^N, 
\label{eq:RTC_tau_scaled}
\end{equation}
where   all other notation is as before.  We again stress that the models \eqref{eq:RTC_tau} and 
\eqref{eq:RTC_tau_scaled} are equivalent, with \eqref{eq:RTC_tau} usually giving the counts of each component and   \eqref{eq:RTC_tau_scaled} providing the normalized abundances.

\begin{remark}
  Historically, the time discretization parameter for the methods described in this paper has
  been $\tau$, leading to the name ``$\tau$-leaping methods.''  We choose to break
  from this tradition
so as not to
  confuse $\tau$ with a stopping time, 
 and we denote our time-step by $h$ 
 to be consistent with much of the numerical analysis literature.
\end{remark}

\section{A review of multi-level Monte Carlo}
\label{sec:old_methods}

Given a stochastic process, $X(\cdot)$, 
let $f: \R^d \to \R$ be a function of the state of the system which gives a measurement of interest.  
Our task is to 
approximate $\E f(X(T))$
efficiently.  
As discussed in Section \ref{sec:intro}, using the ``crude Monte Carlo'' estimator \eqref{eq:crude_approx} with a weakly first order method will provide an estimate with an accuracy of $O(\epsilon)$, in the sense of confidence intervals,  at a computational cost of $O(\epsilon^{-3})$.

%

In multi-level Monte Carlo (MLMC) paths of varying step-sizes are generated and are coupled in
an intelligent manner so that the computational complexity is reduced to 
$O(\epsilon^{-2}(\log \epsilon)^2)$ \cite{Giles2008}. 
Sometimes even the $\log(\epsilon)$ terms can be reduced further \cite{Giles2007}.  
Suppose we have an approximate method, such as Euler's method in the diffusive setting, which is known to be first order accurate in a weak sense, and 1/2 order accurate in a strong $L^2$ sense.  The MLMC estimator  is then built in the following manner.   For a fixed integer $M$, and $\ell \in \{0,1,\dots,L\}$, where $L$ is to be determined, let $\hl = T M^{-\ell}$.  Reasonable choices for $M$ include 2, 3, and 4.  We will denote $Z_{\ell}$ as the approximate process generated using a step-size of $h_{\ell}$.  Choose $L = O(\ln(\epsilon^{-1})),$ so that $h_L = O(\epsilon)$ and   $\E f(X(T)) - \E f(Z_{L}(T)) = O(\epsilon)$, and the bias (i.e. the first term on the right hand side of \eqref{eq:bias_sd}) is of the desired order of magnitude. We  then have
\begin{align}
  \E f(Z_{L}(T)) = \E [ f(Z_{0}(T)) ] + \sum_{\ell = 1}^L \E[ f(Z_{\ell}(T)) - f(Z_{\ell - 1}(T))],
  \label{eq:add_sub}
\end{align}
where the telescoping sum is the key feature to note.
We will now denote the estimator of $\E [ f(Z_{0}(T)) ] $ using $n_{0}$ paths by $\Ql_{0}$, and the estimator of $\E[ f(Z_{\ell}(T)) - f(Z_{\ell - 1}(T))]$ using $n_{\ell}$ paths as $\Ql_{\ell}$.  That is
\begin{align}
	\Ql_{0} &\eqdef \frac{1}{n_{0}} \sum_{i = 1}^{n_{0}} f(Z_{0,[i]}(T)), \quad \text{and} \quad 
	\Ql_{\ell} \eqdef \frac{1}{n_{\ell}} \sum_{i = 1}^{n_{\ell}} ( f(Z_{\ell,[i]}(T)) - f(Z_{\ell - 1,[i]}(T))),
	\label{eq:QL}
\end{align}
where the important point is that both $ Z_{\ell,[i]}(T)$ and $Z_{\ell - 1,[i]}(T)$ are generated using the same randomness, but are constructed using different time discretizations (see \cite{Giles2008, HighamIMA2011} for details on how to do this in the diffusive setting).   We then let 
\begin{equation}
  \Ql \ \eqdef \   \sum_{\ell = 0}^L \Ql_{\ell},
  \label{eq:estimator}
\end{equation}
be the unbiased estimator for $\E[f(Z_{L}(T))]$.  Assuming that we can show 
$\mathsf{Var}(f(Z_{\ell}(T)) - f(Z_{\ell-1}(T))) = O(h_{\ell})$,
 which follows if the method has a strong error of order 1/2 and $f$ is Lipschitz, we may set
\begin{equation*}
	n_{\ell} = O(\epsilon^{-2} Lh_{\ell}),
\end{equation*}
  which yields
$\mathsf{Var}(\Ql) = O(\epsilon^2),$ but with a total computational complexity of $O(\epsilon^{-2}(\log \epsilon)^2)$. 
We make the following observations.  
\begin{enumerate}
\item The gains in computational efficiency come about  for two reasons.  First, a coordinated sequence of simulations are being done, with nested step-sizes, and the simulations with larger step-size are much cheaper than those with very fine step sizes.  Second, while we do still require the generation of paths with fine step-sizes, the variance of $f(Z_{\ell}) - f(Z_{\ell-1})$ will be  small, thereby requiring 
significantly fewer of these expensive paths in the estimation of $\Ql_{\ell}$ of \eqref{eq:QL}.
\item For the analysis in \cite{Giles2008}, it is necessary 
to know both the weak (for the choice of $h_L$) and strong  (for the variance of $\Ql_{\ell}$) behavior of the numerical method, even though we are only solving the {\em weak} approximation problem.
\item The estimator \eqref{eq:estimator} is a biased estimator of $\E f(X(T))$, and the number of levels $L$ was chosen 
to ensure that the bias is within the 
desired tolerance.
\end{enumerate}

\section{Multi-level Monte Carlo for continuous time Markov chains}
\label{sec:new_methods}

We now consider the problem of estimating $\E f(X^N(T))$, where $X^N$ satisfies the general system \eqref{eq:main_multi}.  We again stress that as $X^N$ of \eqref{eq:main_multi} is equivalent to the process $X$ of \eqref{eq:RTC_1st}, efficiently approximating values of the form $\E f(X^N(T))$, for suitable $f$, is \textit{equivalent} to efficiently approximating values of the form $\E g(X(T))$, for suitable functions $g$.  The scaled systems are easier to analyze because the temporal and other quantitative scales have been made explicit.  

Recall that $\N = N^{\gamma} \sum_k N^{c_k}$ gives the order of magnitude of the number of steps needed to generate a single path using an exact algorithm.
As discussed in Section \ref{sec:intro}, to approximate $\E f(X^N(T))$ to an order of accuracy of $\epsilon>0$ using an exact algorithm (such as Gillespie's algorithm or the next reaction method) combined with the crude Monte Carlo estimator, we  need to  generate $\epsilon^{-2}$ paths.  Thus, we have a total computational complexity of $O(\N \epsilon^{-2})$ .   

We will now extend the core ideas of multi-level Monte Carlo as described in Section \ref{sec:old_methods} to the continuous time Markov chain setting with Euler tau-leaping, given in  \eqref{eq:RTC_tau_scaled}, as our approximation method.  We again fix an integer $M>0$, and for $\ell \in \{\ell_0,\dots,L\},$ where both $\ell_0$ and $L$ are to be determined, let $h_{\ell} = TM^{-\ell}$.  We then denote by $Z^N_{\ell}$ the approximate process \eqref{eq:RTC_tau_scaled} generated with a step-size of  $h_{\ell}$.  By \cite{AndMaso2011}, for suitable $f$
\begin{equation*}
  \E f(X^N(T)) - \E f(Z_{\ell}^N(T)) = O(h_{\ell}).
\end{equation*}
   Choose $L = O(\ln(\epsilon^{-1}))$, 
so that $h_L = O(\epsilon)$ and the bias is of the desired 
order of magnitude. We then introduce another telescoping sum 
\begin{align}\label{eq:add_sub_real}
  \E f(Z^N_{L}(T)) = \E [ f(Z^N_{\ell_0}(T)) ] + \sum_{\ell = \ell_0+1}^L \E[ f(Z^N_{\ell}(T)) - f(Z^N_{\ell - 1}(T))].
\end{align}
We will again denote the estimator of 
$\E [ f(Z^N_{\ell_0}(T)) ] $ using $n_{0}$ paths by 
$\Ql_{0}$, and the estimator of 
$\E[ f(Z^N_{\ell}(T)) - f(Z^N_{\ell - 1}(T))]$ 
using $n_{\ell}$ paths by $\Ql_{\ell}$.  
That is
\begin{align}
	\Ql_{0} &\eqdef \frac{1}{n_{0}} \sum_{i = 1}^{n_{0}} f(Z^N_{\ell_0,[i]}(T)), \quad \text{and} \quad 
	\Ql_{\ell} \eqdef \frac{1}{n_{\ell}} \sum_{i = 1}^{n_{\ell}} ( f(Z^N_{\ell,[i]}(T)) - f(Z^N_{\ell - 1,[i]}(T))),
	\label{eq:QL_real}
\end{align}
where we hope that $Z^N_{\ell,[i]}$ and $Z^N_{\ell - 1,[i]}$ can be generated in such a way that $\mathsf{Var}(\Ql_{\ell})$ is small.    We will then let 
\begin{equation}
  \Ql \ \eqdef \   \Ql_{0} +  \sum_{\ell = \ell_0+1}^L \Ql_{\ell},
  \label{eq:estimator_real}
\end{equation}
be the unbiased estimator for $\E[f(Z^N_{L}(T))]$.  The choices for $n_{\ell}$ will depend upon the variances of $\Ql_{\ell}$.

The main requirements for effectively extending MLMC to the current setting now come into focus.  First, we must be able to 
simulate the paths  $Z^N_{\ell}$ and $Z^N_{\ell - 1}$ 
simultaneously in a manner that is efficient and
produces small variances between the paths.   
Second, we must be able to quantify this variance in order 
to control the variance of the associated $\Ql_{\ell}$ terms of \eqref{eq:QL_real}.  Both requirements demand a good coupling of the processes $Z^N_{\ell}$ and $Z^N_{\ell - 1}$.

We motivate our choice of coupling by first treating two simpler tasks.  
First, consider the problem of trying to understand 
the difference between $Z_1(t)$ and $Z_2(t)$, 
where $Z_1,Z_2$ are Poisson processes with rates 
$13.1$ and $13$, respectively. 
A simple approach is to let 
$Y_1$ and $Y_2$ be independent, 
unit-rate Poisson processes,  set 
\[
	Z_1(t) = Y_1(13.1 t) \quad \mathrm{and} \quad 
	Z_2(t) = Y_2(13t),
\]
and consider $Z_1(t) - Z_2(t)$.  
Using this representation, these processes are independent 
and, hence, not coupled.  Further, the variance of their difference is the sum of their variances, and so
\[
\textsf{Var}(Z_1(t) - Z_2(t)) = \textsf{Var}(Z_1(t)) + \textsf{Var}(Z_2(t)) = 26.1t.
\]
Another choice is to let $Y_1$  and $Y_2$ be independent unit-rate Poisson processes, and set 
\[
	Z_1(t) = Y_1(13 t) + Y_2(0.1 t)  
                                 \quad \mathrm{and} \quad
	Z_2(t) = Y_1(13t),
\]
where we have used the additivity property of Poisson 
processes.  
The important point to note is that both $Z_1$ and $Z_2$ 
are using the process $Y_1(13t)$ to generate simultaneous 
jumps.  The process $Z_1$ then uses the auxiliary process 
$Y_2(0.1 t)$ to jump the extra times that $Z_2$ does not.  
The processes $Z_1$ and $Z_2$ will jump together the vast 
majority of times, and hence are tightly coupled; by construction $\textsf{Var}(Z_1(t) - Z_2(t)) = \textsf{Var}(Y_2(0.1 t)) = 0.1t$. 
More generally, if $Z_1$ and $Z_2$ are instead 
inhomogeneous Poisson processes with intensities 
$f(t)$ and $g(t)$, respectively, then we could let 
$Y_1$, $Y_2$, and $Y_3$ be independent, 
unit-rate Poisson processes and define
\begin{align*}
	Z_1(t) &= Y_1\left(\int_0^t f(s) \wedge g(s) ds\right) 
       + Y_2\left( \int_0^t f(s) - 
           \left( f(s) \wedge g(s)  \right) ds\right),\\
	Z_2(t) &=  Y_1\left(\int_0^t f(s) \wedge g(s) ds\right) 
     + Y_3\left( \int_0^t g(s) - \left(f(s) \wedge g(s)  
                           \right)ds\right),
\end{align*}
where we are using that, for example,
\begin{equation*}
Y_1\left(\int_0^t f(s) \wedge g(s) ds\right) + 
        Y_2\left( \int_0^t f(s) - 
             \left( f(s) \wedge g(s)\right) ds\right) 
            \eqdist Y\left( \int_0^t f(s) ds\right),
\end{equation*}
where $Y$ is a unit rate Poisson process and we recall that 
$a\wedge b \eqdef \min\{a,b\}$.

We now return to the main problem of coupling the processes 
$Z_{\ell}^N$ and $Z_{\ell-1}^N$, each satisfying 
\eqref{eq:RTC_tau_scaled} with their respective step-sizes.   
We couple the processes $Z^N_{\ell}$ and $Z^N_{\ell-1}$ 
in the following manner, which is similar to a 
coupling originally  used in \cite{Kurtz82}, 
and later in  \cite{AndersonGangulyKurtz}, as an analytical 
tool, and subsequently in \cite{AndersonSensitivities} towards the problem of computing parameter sensitivities:
\begin{align}
\begin{split}
	\Zm_{\ell}(t)  &= \Zm_{\ell}(0) + \sum_{k}  Y_{k,1}\left(N^{\gamma}N^{c_k}  \int_0^t \lambda_k(\Zm_{\ell} \circ \eta_{\ell}(s)) \wedge \lambda_k(\Zm_{\ell-1} \circ \eta_{\ell-1}(s)) ds \right)\zeta_k^N \\
	& + \sum_{k} Y_{k,2}\left(N^{\gamma}N^{c_k} \int_0^t  \lambda_k(\Zm_{\ell}\circ \eta_{\ell}(s)) -  \lambda_k(\Zm_{\ell} \circ \eta_{\ell}(s)) \wedge \lambda_k(\Zm_{\ell-1} \circ \eta_{\ell-1}(s)) ds \right)\zeta_k^N,  \label{eq:Z1}
    \end{split}\\
    \begin{split}
	\Zm_{\ell-1}(t) &= \Zm_{\ell-1}(0) + \sum_{k} Y_{k,1}\left( N^{\gamma}N^{c_k}  \int_0^t \lambda_k(\Zm_{\ell} \circ \eta_{\ell}(s)) \wedge \lambda_k(\Zm_{\ell-1} \circ \eta_{\ell-1}(s))  ds  \right)\zeta_k^N \\
	& + \sum_k Y_{k,3}\left( N^{\gamma}N^{c_k} \int_0^t  \lambda_k(\Zm_{\ell-1}\circ \eta_{\ell-1}(s)) -  \lambda_k(\Zm_{\ell} \circ \eta_{\ell}(s)) \wedge \lambda_k(\Zm_{\ell-1} \circ \eta_{\ell-1}(s)) ds \right)\zeta_k^N,
    \label{eq:Z2}
    \end{split}	
\end{align}
 where the $Y_{k,i},\ i \in \{1,2,3\}$, are independent, 
unit-rate Poisson processes, and for each $\ell$, 
we define 
$\eta_{\ell}(s) \eqdef \lfloor s/h_{\ell}\rfloor h_{\ell}$.   
Note that we essentially used the coupling of the simpler 
examples above (pertaining to $Z_1$ and $Z_2$) 
for each of the reaction channels.   
 
 The paths of the coupled processes can easily be computed simultaneously and the  distributions of the marginal processes are the same as the usual scaled Euler approximate paths 
\eqref{eq:RTC_tau_scaled} with similar step-sizes.  
More precisely, the system \eqref{eq:Z1}--\eqref{eq:Z2} is the 
scaled version of, and is hence equivalent to, the system
 \begin{align}
 \begin{split}
	Z_{\ell}(t) = &Z_{\ell}(0) + \sum_{k} Y_{k,1}\left( \int_0^t \lambda_k(Z_{\ell}\circ \eta_{\ell}(s)) \wedge \lambda_k(Z_{\ell-1}\circ \eta_{\ell-1}(s)) ds \right)\zeta_k \\
	&  + \sum_{k} Y_{k,2}\left(  \int_0^t  \lambda_k(Z_{\ell}\circ \eta_{\ell}(s)) - \lambda_k(Z_{\ell}\circ \eta_{\ell}(s)) \wedge \lambda_k(Z_{\ell-1}\circ \eta_{\ell-1}(s)) ds \right)\zeta_k,\\
	Z_{\ell-1}(t) = &Z_{\ell-1}(0) + \sum_{k} Y_{k,1}\left( \int_0^t \lambda_k(Z_{\ell}\circ \eta_{\ell}(s)) \wedge \lambda_k(Z_{\ell-1}\circ \eta_{\ell-1}(s))   ds  \right)\zeta_k\\
	&  + \sum_{k} Y_{k,3}\left(\int_0^t \lambda_k( Z_{\ell-1}\circ \eta_{\ell-1}(s)) - \lambda_k(Z_{\ell}\circ \eta_{\ell}(s)) \wedge \lambda_k(Z_{\ell-1}\circ \eta_{\ell-1}(s)) ds \right)\zeta_k,
	\end{split}
	\label{eq:unscaled_coupled_tau}
\end{align}
where now the marginal processes are distributionally equivalent to the approximate processes \eqref{eq:RTC_tau} with similar step-sizes, and all notation is as before.
The natural algorithm to simulate the representation \eqref{eq:unscaled_coupled_tau} (and hence \eqref{eq:Z1}--\eqref{eq:Z2}) 
to a time $T>0$ is the following.

\begin{algorithm}[Simulation of the representation \eqref{eq:unscaled_coupled_tau}]
Fix an integer $M\ge 2$.  Fix $h_{\ell}>0$ and set $h_{\ell-1} = M\times h_{\ell}$.  Set $Z_{\ell}(0) = Z_{\ell-1}(0) = x_0$, $t_0 = 0$,
  $n=0$.  Repeat the following steps until $t_{n} \ge T$:
  
  \begin{enumerate}[$(i)$]
  \item For $j = 0,\dots,M-1$,

 	 \begin{enumerate}
  		\item Set 
  		\begin{itemize}
  			\item $A_{k,1} = \lambda_k(Z_{\ell}( t_n + j\times h_{\ell})) \wedge \lambda_k(Z_{\ell-1}(t_n))$.
			\item $A_{k,2} =  \lambda_k(Z_{\ell}(t_n + j\times h_{\ell})) - A_{k,1}$.
			\item $A_{k,3} = \lambda_k(Z_{\ell-1}(t_n)) - A_{k,1}$.
  		\end{itemize}
		\item  For each $k$, let 
		\begin{itemize}
		\item $\Lambda_{k,1} = \text{Poisson}(A_{k,1} h_{\ell})$.
		\item $\Lambda_{k,2} = \text{Poisson}(A_{k,2} h_{\ell})$.
		\item$\Lambda_{k,3} = \text{Poisson}(A_{k,3} h_{\ell})$.
		\end{itemize}
		
		\item Set 
		\begin{itemize}
			\item $Z_{\ell}(t_n + (j+1)\times h_{\ell}) = Z_{\ell}(t_n + j\times h_{\ell}) + \sum_k (\Lambda_{k,1} + \Lambda_{k,2}) \zeta_k.$
			\item  $Z_{\ell-1}(t_n + (j+1)\times h_{\ell}) = Z_{\ell-1}(t_n + j\times h_{\ell}) + \sum_k (\Lambda_{k,1} + \Lambda_{k,3}) \zeta_k.$
		\end{itemize}
  	\end{enumerate}
   \item Set $t_{n+1} = t_{n} + h_{\ell-1}$.  
   \item Set $n \leftarrow n+1$.
  \end{enumerate}
  \label{alg:double_tau}
\end{algorithm}

We make the following observations.  
First, while Algorithm \ref{alg:double_tau} formally simulates 
the representation \eqref{eq:unscaled_coupled_tau}, 
the scaled version of the process generated via 
Algorithm \ref{alg:double_tau} satisfies 
\eqref{eq:Z1}--\eqref{eq:Z2}.  
Second, we do not need to update either 
$Z_{\ell-1}$ or $\lambda_k(Z_{\ell-1})$ during the workings 
of the inner loop of $j = 0,\dots, M-1$.  
Third, at most one of $A_2,A_3$ will be non-zero during 
each step, with both being zero whenever 
$\lambda_k(Z_{\ell}(t_n)) = \lambda_k(Z_{\ell-1}(t_n))$.  
Therefore, at most two Poisson random variables will be 
required per reaction channel at each step and not three.   
Fourth, the above algorithm, and hence the couplings 
\eqref{eq:unscaled_coupled_tau} and/or 
\eqref{eq:Z1}--\eqref{eq:Z2}, is no harder to simulate, 
from an implementation standpoint, than the usual 
Euler tau-leaping.  
Fifth, while two paths are being generated, it should be the 
case that $\max\{A_2,A_3\}$ is small for each step.  
Hence the work in computing the Poisson random variables 
will fall on $\Lambda_{k,1}$,\footnote{The cost of generating a Poisson random variable generally increases with the size of the mean} which is the same amount of 
work as would be needed for the generation of a 
\textit{single} path of Euler tau-leaping.  
 
 In Section \ref{sec:analysis} we will prove the following theorem, which is one of our main analytical results.
 \begin{theorem} \label{thm:MLMC}
    	Suppose $(\Zm_{\ell},\Zm_{\ell-1})$ satisfy \eqref{eq:Z1} and \eqref{eq:Z2} with $\Zm_{\ell}(0)=\Zm_{\ell-1}(0)$.  Then, there exist functions $C_1,C_2$, that do not depend on $h_{\ell}$, such that 
   \begin{equation*}
      \sup_{t \le T} \E| \Zm_{\ell}(t) - \Zm_{\ell - 1}(t) |^2 \le C_1(N^{\gamma}T)N^{-\rho}(N^{\gamma} h_{\ell}) + C_2(N^{\gamma}T)   (N^{\gamma} h_{\ell})^2.
         \end{equation*} 
In particular, for $\gamma \le 0$ the values 
$
C_1(N^{\gamma}T)
$ 
and
$
C_2(N^{\gamma}T)
$ 
may be bounded above uniformly in $N$.
 \end{theorem}

\begin{remark}\label{remark:Cs}
 The specific forms of $C_1(N^{\gamma}T)$ and $C_2(N^{\gamma}T)$ 
for  Theorem \ref{thm:MLMC} and Theorem \ref{thm:MLMC_exact} 
below are given in Section \ref{sec:analysis}.  
However, we note here that if $\gamma>0$, the factors 
$C_1(N^{\gamma}T)$ and $C_2(N^{\gamma}T)$ could be huge, leading to
 upper bounds in Theorem~\ref{thm:MLMC}
and Theorem~\ref{thm:MLMC_exact} of no practical use.   
So we henceforth assume that
$\gamma \le 0$ and thus regard 
$C_1(N^{\gamma}T)$ and $C_2(N^{\gamma}T)$
as constants independent of $N$.
We note that the classical chemical kinetics scaling, with $\gamma = 0$, 
satisfies this assumption. 
However,
 good performance is observed 
 in Section \ref{sec:examples}
 with 
 $\gamma>0$, suggesting that 
 further analysis may extend the
 range of validity for this method.
 \end{remark}
 
 Note that Theorem \ref{thm:MLMC} together with $f$ Lipschitz gives us the estimate
 \begin{eqnarray}
 \left|  
 \mathsf{Var}\left( f(  \Zm_{\ell}(t)  ) \right)
 - 
  \mathsf{Var}\left( f(  \Zm_{\ell-1}(t)  ) \right)
\right|
 &\le&
  \E 
   \left|
     f(  \Zm_{\ell}(t)  )
  - 
   f(  \Zm_{\ell-1}(t)  )
   \right|^2
  \nonumber \\
 &\le&  
 C
     \E| \Zm_{\ell}(t) - \Zm_{\ell - 1}(t) |^2 
 \nonumber \\
   &\le&  
  C \left[ C_1(N^{\gamma}T)N^{-\rho}(N^{\gamma}h_{\ell}) + C_2(N^{\gamma}T)   (N^{\gamma}h_{\ell})^2 \right],
\label{eq:varest}
 \end{eqnarray}
 which we will use to control the variance of $\Ql_{\ell}$ in \eqref{eq:QL_real}.
 
Before further exploring MLMC in the current setting, we present a coupling of the exact process $X^N$ and the approximate process $Z^N_{\ell}$.  We will later use this coupling to produce an unbiased MLMC estimator.    We define $X^N$ and $Z^N_{\ell}$ via
\begin{align}
\begin{split}
	\X(t) = &X^N(0) + \sum_{k} Y_{k,1}\left( N^{\gamma}N^{c_k} \int_0^t  \lambda_k(\X(s))  \wedge  \lambda_k(\Zm_{\ell}\circ \eta_{\ell}(s))  ds  \right)\zeta_k^N\\
	&\hspace{.1in} + \sum_{k} Y_{k,2}\left( N^{\gamma}N^{c_k} \int_0^t \lambda_k(\X(s))  -  \lambda_k(\X(s))  \wedge   \lambda_k(\Zm_{\ell}\circ \eta_{\ell}(s))  ds \right)\zeta_k^N,
    \label{eq:Z_X1}
\end{split}\\
\begin{split}
	\Zm_{\ell}(t) = &\Zm_{\ell}(0) + \sum_{k} Y_{k,1}\left(N^{\gamma}N^{c_k}  \int_0^t \lambda_k(\X(s))  \wedge   \lambda_k(\Zm_{\ell}\circ \eta_{\ell}(s))  ds \right)\zeta_k^N \\
	&\hspace{.1in} + \sum_{k} Y_{k,3}\left(N^{\gamma}N^{c_k}  \int_0^t  \lambda_k(\Zm_{\ell}\circ \eta_{\ell}(s))  - \lambda_k(\X(s))  \wedge  \lambda_k(\Zm_{\ell}\circ \eta_{\ell}(s)) ds  \right)\zeta_k^N, 
      \label{eq:Z_X2}
    \end{split}	
\end{align}
where all notation is as before.  Note that the distributions of the marginal processes $X^N$ and $Z^N_{\ell}$  are equal to those of  \eqref{eq:main_multi} and \eqref{eq:RTC_tau_scaled}.  The unscaled processes satisfy
\begin{align}
\begin{split}
	X(t) = &X(0) + \sum_{k} Y_{k,1}\left( \int_0^t  \lambda_k(X(s))  \wedge  \lambda_k(Z_{\ell}\circ \eta_{\ell}(s))  ds  \right)\zeta_k\\
	&\hspace{.1in} + \sum_{k} Y_{k,2}\left(  \int_0^t \lambda_k(X(s))  -  \lambda_k(X(s))  \wedge   \lambda_k(Z_{\ell}\circ \eta_{\ell}(s))  ds \right)\zeta_k,\\
	Z_{\ell}(t) = &Z_{\ell}(0) + \sum_{k} Y_{k,1}\left(  \int_0^t \lambda_k(X(s))  \wedge   \lambda_k(Z_{\ell}\circ \eta_{\ell}(s))  ds \right)\zeta_k \\
	&\hspace{.1in} + \sum_{k} Y_{k,3}\left(\int_0^t  \lambda_k(Z_{\ell}\circ \eta_{\ell}(s))  - \lambda_k(X(s))  \wedge  \lambda_k(Z_{\ell}\circ \eta_{\ell}(s)) ds  \right)\zeta_k, 
    \end{split}
    \label{eq:unscaled_coupled_exact}
\end{align}
which is equivalent to \eqref{eq:Z_X1} and \eqref{eq:Z_X2}, and whose marginal processes have the same distributions as \eqref{eq:RTC_1st} and \eqref{eq:RTC_tau}. 

 The natural algorithm to simulate \eqref{eq:unscaled_coupled_exact}, and hence
 \eqref{eq:Z_X1}--\eqref{eq:Z_X2}, is the next reaction method 
\cite{Anderson2007a,Gibson2000}, where the system 
is viewed as having 
dimension
 $2 d$ 
with state $(X^N,Z_{\ell}^N)$, and each of the ``next reactions'' must be calculated over the Poisson processes $Y_{k,1},Y_{k,2},Y_{k,3}$.  See \cite{Anderson2007a} for a thorough explanation of how the next reaction method is equivalent to simulating representations of the forms considered here.    Below, we will denote a uniform$[0,1]$ random variable by rand$(0,1)$, and we remind the reader that if $U \sim \text{rand}(0,1)$, then $\ln(1/U)$ is an exponential random variable with a parameter of one.  All random variables generated are assumed to be independent of each other and all previous random variables.

\begin{algorithm}
   [Simulation of the representation \eqref{eq:unscaled_coupled_exact}]
\textbf{Initialize}.  Fix $h_{\ell}>0$.  Set $X(0) = Z_{\ell}(0)  =x_0$ and $t = 0$.  Set $\widetilde Z_{\ell} = Z_{\ell}(0)$.  Set $T_{\mathrm{tau}} = h_{\ell}$.   For each $k\in \{1,\dots,R\}$ and $i \in \{1,2,3\}$, set $P_{k,i} = \ln(1/r_{k,i})$, where $r_{k,i}$ is rand$(0,1)$, and $T_{k,i} = 0$.
  \begin{enumerate}[$(i)$]
  \item For each $k$, set
  	\begin{itemize}
  		\item $A_{k,1} = \lambda_k(X(t)) \wedge \lambda_k(\widetilde Z_{\ell})$.
		\item $A_{k,2} =  \lambda_k(X(t)) - A_{k,1}$.
		\item $A_{k,3} = \lambda_k(\widetilde Z_{\ell}) - A_{k,1}$.
  	\end{itemize}
\item  For each $k\in \{1,\dots,R\}$ and $i \in \{1,2,3\}$, set 
\begin{equation*}
   \Delta t_{k,i} = \left\{ \begin{array}{cr}
   (P_{k,i} - T_{k,i})/A_{k,i}, & \text{ if } A_{k,i} \ne 0\\
   \infty, & \text{ if } A_{k,i} = 0
   \end{array} \right. .
 \end{equation*}
\item Set $\Delta = \min_{k,i}\{\Delta t_{k,i}\}$, and let $\mu \equiv \{k,i\}$ be the indices where the minimum is achieved.
	\item If $t + \Delta \ge  T_{\mathrm{tau}}$, 
	\begin{enumerate}
		\item Set $\widetilde Z_{\ell} =  Z_{\ell}(t)$.
		\item For each $k\in \{1,\dots,R\}$ and $i \in \{1,2,3\}$, set $T_{k,i} = T_{k,i} + A_{k,i}\times (T_{\mathrm{tau}} - t)$.
		\item Set $t = T_{\mathrm{tau}}$.
		\item Set $T_{\mathrm{tau}} = T_{\mathrm{tau}} + h_{\ell}$.
		\item Return to step $(i)$ or quit.
	\end{enumerate} 
	\item Else, 
	\begin{enumerate}
		\item Update.  For $\{k,i\} = \mu$, where $\mu$ is from $(iii)$,
		\begin{itemize}
		\item If $i = 1$, set $X(t+\Delta) = X(t) + \zeta_k$ and $Z_{\ell}(t + \Delta) = Z_{\ell}(t) + \zeta_k$.
		\item If $i = 2$, set $X(t+\Delta) = X(t) + \zeta_k$.
		\item If $i = 3$, set $Z_{\ell}(t + \Delta) = Z_{\ell}(t) + \zeta_k$
		\end{itemize}
		\item For each $k\in \{1,\dots,R\}$ and $i \in \{1,2,3\}$, set $T_{k,i} = T_{k,i} + A_{k,i}\times \Delta$.
		\item Set $P_{\mu} = P_{\mu} + \ln(1/r)$, where $r$ is rand$(0,1)$, and $\mu$ is from $(iii)$.
		\item Set $t = t + \Delta$.
		\item Return to step $(i)$ or quit.
	\end{enumerate}
  \end{enumerate}
\end{algorithm}

The following theorem, which should be compared with Theorem \ref{thm:MLMC}, is proven in Section \ref{sec:analysis} and is our second main analytical result.
\begin{theorem} \label{thm:MLMC_exact}
    	Suppose $(X^N,\Zm_{\ell})$ satisfy \eqref{eq:Z_X1} and \eqref{eq:Z_X2} with $X^N(0) = \Zm_{\ell}(0)$.  Then, there exist functions $C_1,C_2$, that 
do not depend on $h_{\ell}$, such that 
   \begin{equation*}
      \sup_{t \le T} \E| X^N(t) - Z^N_{\ell}(t) |^2 \le C_1(N^{\gamma}T)N^{-\rho} (N^{\gamma} h_{\ell})
 + C_2(N^{\gamma}T)  (N^{\gamma} h_{\ell})^2.
\end{equation*} 
Moreover, for $\gamma \le 0$ the values 
$
C_1(N^{\gamma}T)
$ 
and
$
C_2(N^{\gamma}T)
$ 
may be bounded above uniformly in $N$.         
 \end{theorem}

 We are now in a position to develop MLMC in the stochastic chemical kinetic setting.  Recall our assumption that $\gamma \le 0$, so $C_1$ and $C_2$ in Threorems~\ref{thm:MLMC} and
      \ref{thm:MLMC_exact} are bounded.  We return to the $\Ql_{\ell}$ terms in  \eqref{eq:QL_real}.  Supposing that the test function $f$ is uniformly Lipschitz in our domain of interest (note that this is automatic for any reasonable $f$ in the case when mass is conserved), then for $\ell > \ell_0$, we know from (\ref{eq:varest}) that 
 \[
 	\mathsf{Var}(\Ql_{\ell}) 
      \le C \frac{1}{n_{\ell}} \left[C_1(N^{\gamma}T)N^{-\rho}(N^{\gamma} h_{\ell}) + C_2(N^{\gamma}T)   (N^{\gamma} h_{\ell})^2\right].
 \] 
Note that if $N^{-\rho} \le h_{\ell}$, the leading order of the error is the $h_{\ell}^2$ term.  As a heuristic argument for this behavior, note that if $N^{-\rho} \le h_{\ell}$ and $N$ is large while $h_{\ell}$ is small,  then the processes are nearing a scaling regime in which deterministic dynamics would be a good approximation for the model $X^N$.  In this case, one should expect that the \textit{squared} difference between two Euler paths should behave like the usual order one error, squared.
 
 We may now conclude that the variance of the estimator $\Ql$ defined in \eqref{eq:estimator_real} satisfies
 \begin{align*}
    \mathsf{Var}(\Ql) &= \mathsf{Var}(\Ql_{\ell_0}) + \sum_{\ell = \ell_0+1}^L \mathsf{Var} (\Ql_{\ell})\\
    &\le \frac{K_0}{n_0} + \sum_{\ell = \ell_0 + 1}^L C \frac{1}{n_{\ell}} \left[C_1(N^{\gamma}T)N^{-\rho}(N^{\gamma} h_{\ell}) + C_2(N^{\gamma}T)   (N^{\gamma} h_{\ell})^2\right],
 \end{align*}
 where $K_0 = \mathsf{Var}(f(Z^N_{\ell_0}(T)))$.
For $h>0$ we define
\begin{equation}
	A(h) \eqdef N^{-\rho}(N^{\gamma}h) + (N^{\gamma} h)^2.
	\label{def:A}
\end{equation}
Letting $n_0 = O(\epsilon^{-2})$, and for $\ell > \ell_0$ letting
\begin{equation*}
	n_{\ell} = O\left(\epsilon^{-2}(L - \ell_0)A(h_{\ell})\right),
\end{equation*}
we see that
\begin{equation*}
  \mathsf{ Var}(\Ql) = O(\epsilon^{2}).
\end{equation*}
 As the computational complexity of generating a single path of the coupled processes $(Z^N_{\ell},Z^N_{\ell-1})$ is $O(h_{\ell}^{-1})$, we see that the total computational complexity of the method with these choices of $n_{\ell}$ is of order
 \begin{align}
 	n_0h_{\ell_0}^{-1} + \sum_{\ell = \ell_0 + 1}^L n_{\ell} h_{\ell}^{-1} &= \epsilon^{-2}h_{\ell_0}^{-1} + \sum_{\ell = \ell_0 + 1}^L\epsilon^{-2}(L - \ell_0)A(h_{\ell})  h_{\ell}^{-1} \notag \\
	&= \epsilon^{-2}\left( h^{-1}_{\ell_0} + (L - \ell_0) \sum_{\ell = \ell_0+1}^L (N^{-\rho}N^{\gamma} + h_{\ell}N^{2\gamma})  \right)\notag \\
	&\le \epsilon^{-2}\left( h^{-1}_{\ell_0} + \ln(\epsilon)^2 N^{-\rho}N^{\gamma} + \ln(\epsilon^{-1})\frac{1}{M-1}h_{\ell_0 }N^{2\gamma}\right)
	\label{eq:CC_1},
 \end{align}
where we used that 
\begin{equation}
	\sum_{\ell = \ell_0 + 1}^L h_{\ell} \le h_{\ell_0} \sum_{\ell = 1}^{\infty} \frac{1}{M^{\ell}}  =h_{\ell_0}  \frac{1}{M-1}.
	\label{eq:geometric_bound}
\end{equation}
A more careful choice of $n_{\ell}$ can potentially reduce the $\ln(\epsilon)$ terms further, see for example \cite{Giles2008}, but in the present case, the computational complexity will be dominated by $\epsilon^{-2}h_{\ell_0}^{-1}$ in most nontrivial examples.  Further, as will be 
 discussed in Section~\ref{sec:implementation}, the $n_{\ell}$ can be 
chosen algorithmically, by optimizing for a given problem.

\subsection{An unbiased MLMC}
\label{sec:unbiased}

We now build an unbiased MLMC estimator for $\E f(X^N(T))$ in a similar manner as before with a single important difference: at the finest scale, we couple $X^N$ with $Z^N_{L}$.  That is, we 
use the identity 
\begin{equation*}
	\E f(X^N(T)) = \E [f(X^N(T)) - f(Z^N_L(T))] + \sum_{\ell = \ell_0+1}^L \E [ f(Z^N_{\ell}) - f(Z^N_{\ell-1})] + \E f(Z^N_{\ell_0}(T)).
\end{equation*}  
For appropriate choices of $n_0,n_{\ell}$, and $n_{E}$, we define the estimators for the three terms above via
\begin{align}
	\Ql_{E} &\eqdef \frac{1}{n_{E}} \sum_{i = 1}^{n_{E}} (f(X^N_{[i]}(T) - f(Z^N_{L,[i]}(T))) ,\notag\\
	\Ql_{\ell} &\eqdef \frac{1}{n_{\ell}} \sum_{i = 1}^{n_{\ell}} ( f(Z^N_{\ell,[i]}(T)) - f(Z^N_{\ell - 1,[i]}(T))), \quad \text{for } \ell \in \{\ell_0+1,\dots,L\},\notag \\
	\Ql_{0} &\eqdef \frac{1}{n_{0}} \sum_{i = 1}^{n_{0}} f(Z^N_{\ell_0,[i]}(T))\notag,
\end{align}
and note that 
\begin{equation}
	\Ql  \eqdef  \Ql_E +  \sum_{\ell = \ell_0+1}^L \Ql_{\ell} + \Ql_0
	\label{eq:unbiased_MLMC}
\end{equation}
is an \textit{unbiased} estimator for $\E f(X^N(T))$.  Applying both Theorems \ref{thm:MLMC} and \ref{thm:MLMC_exact} yields
\begin{align*}
	\mathsf{Var}(\Ql_{E}) &\le K_1(N^{\gamma} T)\frac{1}{n_E} A(h_L),\\
	\mathsf{Var}(\Ql_{\ell}) &\le K_2(N^{\gamma} T) \frac{1}{n_{\ell}}A(h_{\ell}), \quad \text{for } \ell \in \{\ell_0+1,\dots, L\},
 \end{align*} 
 where $K_1(N^{\gamma} T)$ and $K_2(N^{\gamma} T)$ are independent of $h_{\ell}$, and,
under our assumption that $\gamma \le 0$, can be bounded uniformly in $N$.  
  It follows that the choice 
  \begin{align}
 \begin{split}
 	n_E &= O(\epsilon^{-2}A(h_L)),\\
	n_{\ell} &= O\left(\epsilon^{-2}(L - \ell_0)A(h_{\ell})\right), \quad \text{ for } \ell \in \{\ell_0,\dots,L\},\\
	n_0 &= O(\epsilon^{-2}),
	\end{split}
	\label{eq:predicted_nl}
 \end{align}
 gives us
 \begin{align*}
 	\mathsf{Var}(\Ql) &= \mathsf{Var}(\Ql_E) + \sum_{\ell = \ell_0+1}^L \mathsf{Var}(\Ql_{\ell}) + \mathsf{Var}(\Ql_0) \\ 
	&= O(\epsilon^2) + \sum_{\ell = \ell_0 + 1}^L O(\epsilon^2 (L-\ell_0)^{-1}) + O(\epsilon^2)\\
	&= O(\epsilon^2).
 \end{align*}
The computational complexity is now of order
 \begin{align}
 	n_E \overline N +  \sum_{\ell = \ell_0 + 1}^L n_{\ell} &h_{\ell}^{-1} + n_0h_{\ell_0}^{-1}  =  \overline N \epsilon^{-2}A(h_L) + \sum_{\ell = \ell_0 + 1}^L\epsilon^{-2}(L - \ell_0)A(h_{\ell})  h_{\ell}^{-1} + \epsilon^{-2}h_{\ell_0}^{-1}  \notag \\
	&= \epsilon^{-2}\left( \overline N A(h_L) +  (L - \ell_0)\sum_{\ell = \ell_0}^L (N^{-\rho}N^{\gamma} + h_{\ell}N^{2\gamma})   + h^{-1}_{\ell_0} \right)\notag \\
	&\le  \epsilon^{-2}\left( \overline N A(h_L) +h^{-1}_{\ell_0}+  \ln(\epsilon)^2 N^{-\rho}N^{\gamma} + \ln(\epsilon^{-1}) \frac{1}{M-1}h_{\ell_0}N^{2\gamma}   \right),\label{eq:CC_2} 
 \end{align}
where we again made use of the inequality \eqref{eq:geometric_bound}.

\subsection{Some observations}
\label{subsec:obs}

 A few observations are in order.  First, in the above analysis of the unbiased MLMC estimator, the weak error of the process $Z_h^N$ plays \textit{no role}.  Thus, there is no reason to choose $h_{L} = O(\epsilon)$ for a desired accuracy of $\epsilon>0$.  Without having to worry about the bias, we have the opportunity to simply choose $h_{L}$ ``small enough'' for $\mathsf{Var}(X^N(\cdot) - Z^N_{L}(\cdot))$ to be small, which can be approximated with a few preliminary simulations before the full MLMC is carried out (see Section \ref{sec:implementation} for more implementation details).

 Second, one of the main impediments to the use of tau-leaping methods has been the possibility for paths to leave the non-positive orthant.  In fact, there have been multiple papers written on the subject of how to enforce non-negativity of species numbers with \cite{Anderson2007b, Cao2005, Chatterjee2005, Tian2004} representing just a sample.  We note that for the unbiased MLMC estimator \eqref{eq:unbiased_MLMC} it almost does not matter how, or even if, non-negativity is enforced.  So long as the processes are well defined on all of $\Z^d$, for example by defining the intensity functions $\lambda_k$ in some reasonable way, and so long as we can still quantify the relations given in Theorems \ref{thm:MLMC} and \ref{thm:MLMC_exact}, everything above still holds.  The cost, to the user, of poorly defining what happens if $Z_h$ leaves the positive orthant will simply be the need for the generation of more paths to reduce the variance of the (still unbiased) estimator.  Of course, this cost could be quite high as negativity of population numbers can lead to instability if they have not defined the intensity functions outside the positive orthant in a reasonable manner.  However, in Section \ref{sec:implementation} we discuss how intelligent implementation of the method can greatly reduce this cost by ensuring that the approximate paths remain stable with high probability.

 Third, inspecting \eqref{eq:CC_1} and \eqref{eq:CC_2} shows 
that the unbiased MLMC estimator \eqref{eq:unbiased_MLMC} 
has an additional term 
of $O(\overline N A(h_{L}) \epsilon^{-2})$ 
in its computational complexity bound, 
as compared with the biased MLMC estimator \eqref{eq:estimator_real}.  The authors feel that $\overline N A(h_{L})$ would have to be quite substantial to warrant not using the unbiased version.
 
 Fourth, note that we always have the following:
\begin{align}
\begin{split}
	\text{Computational complexity of unbiased MLMC} &= O\left(\epsilon^{-2} (\overline N A(h_L) + h_{\ell_0}^{-1} + \text{log term})\right)\\
	& \ll O\left( \epsilon^{-2} \overline N\right)\\
	&= \begin{array}{l}
	   \text{Computational complexity of exact algorithm}\\
	   \text{with crude Monte Carlo.}
	   \end{array}
	   \end{split}
	   \label{eq:CC}
\end{align}
 Thus, under our standing assumption $\gamma \le 0$, the unbiased MLMC estimator should 
 be the method of choice over using an exact algorithm alone together 
 with crude Monte Carlo, which is by far the most popular method today.
  For example, consider the case when the system satisfies the 
 classical scaling, for which $\rho = 1$, $\gamma = 0$ and 
 $c_k \equiv 1$. 
 In this case, $\overline N = N$ and, as there is little reason to use an approximate method with a time step that is smaller than the order of magnitude of the wait time between jumps for an exact method, we may assume that $h_{L} > 1/N = N^{-\rho}$. Therefore, in this specific case, $A(h_{L}) = O(h_L^2)$ and the computational speedup predicted by \eqref{eq:CC_2} and/or \eqref{eq:CC} is of the order
 \begin{align*}
 	\text{Speed-up factor} \approx \frac{\epsilon^{-2}N}{\epsilon^{-2}(N h_L^2 + h_{\ell_0}^{-1} + \log(\epsilon))} = \frac{N}{N h_L^2 + h_{\ell_0}^{-1} + \log(\epsilon)}. 
 \end{align*}
 Thus we have
 \begin{align*}
 	\text{Speed-up factor}\gtrapprox \min\left( h_L^{-2}, N h_{\ell_0} \right).
 \end{align*}
Therefore, even though the method is unbiased, the computational burden has been shifted from the exact process to that of an approximate process with a crude time-step.  This behavior is demonstrated in an example found in Section \ref{sec:examples}, though on a system not satisfying the classical scaling.
 
 Note also that \eqref{eq:CC} holds even if $\overline N$, the approximate cost of computing a single path, is not extremely large.  For example, even if the cost is only in the hundreds, or maybe thousands, of steps per exact path, the above analysis points out that if great accuracy is required (so that $\epsilon^{-2}$ is very large), the unbiased MLMC estimator will still decrease the computational complexity substantially.  It should be pointed out that in these cases of moderate $\overline N$, we will typically have $\gamma \le 0$ and so the analysis will hold.    

The conclusion of this analysis, 
backed up by the examples in Section \ref{sec:examples}, 
is that MLMC methods, with processes coupled via the 
representations \eqref{eq:unscaled_coupled_tau} 
and \eqref{eq:unscaled_coupled_exact}, 
and the unbiased MLMC in particular,  
produce substantial gains in computational efficiency and 
could become standard algorithms in the sciences.  
Further attention, however, needs to be given to the 
case $\gamma>0$, and this will be
a focus for future work.

\section{Delayed proofs of Theorems \ref{thm:MLMC} and \ref{thm:MLMC_exact}}
\label{sec:analysis}

 We begin by focussing on the proof of Theorem \ref{thm:MLMC_exact}, which is restated here for completeness.\vspace{.1in}
 
\noindent \textbf{Theorem 2.} \textit{
    	Suppose $(X^N,\Zm_{\ell})$ satisfy \eqref{eq:Z_X1} and \eqref{eq:Z_X2} with $X^N(0) = \Zm_{\ell}(0)$.  Then, there exist functions $C_1,C_2$, that 
do not depend on $h_{\ell}$, such that 
   \begin{equation*}
      \sup_{t \le T} \E| X^N(t) - Z^N_{\ell}(t) |^2 \le C_1(N^{\gamma}T)N^{-\rho} (N^{\gamma} h_{\ell})
 + C_2(N^{\gamma}T)  (N^{\gamma} h_{\ell})^2.
\end{equation*} 
Moreover, for $\gamma \le 0$ the values 
$
C_1(N^{\gamma}T)
$ 
and
$
C_2(N^{\gamma}T)
$ 
may be bounded above uniformly in $N$.         
         }

 We start with the following lemma. 

   \begin{lemma}    \label{lem:first_bound_exact}
   Suppose $(X^N,\Zm_{\ell})$ satisfy \eqref{eq:Z_X1} and \eqref{eq:Z_X2} with $X^N(0) = \Zm_{\ell}(0)$.  Then, there exist positive constants $c_1,c_2$, independent of $N$, $\gamma$, and $T$, such that for $t \ge 0$
   \begin{equation*}
      \E| X^N(t) - Z^N_{\ell}(t) | \le c_1  \left(e^{c_2N^{\gamma}t} -1\right) (N^{\gamma} h_{\ell}).
   \end{equation*}
    \end{lemma}
   
   \begin{proof}
   Note that
	\begin{align*}
	 \E |X^N(t) - &Z^N_{\ell}(t)| \\
	 &=  \E \bigg| \sum_{k} Y_{k,2}\left( N^{\gamma}N^{c_k} \int_0^t \lambda_k(\X(s))  -  \lambda_k(\X(s))  \wedge   \lambda_k(\Zm_{\ell}\circ \eta_{\ell}(s))  ds \right) \zeta_k^N   \\
&\hspace{.2in} - \sum_{k} Y_{k,3}\left(N^{\gamma}N^{c_k} \int_0^t  \lambda_k(\Zm_{\ell}\circ \eta_{\ell}(s))  - \lambda_k(\X(s))  \wedge  \lambda_k(\Zm_{\ell}\circ \eta_{\ell}(s)) ds  \right) \zeta_k^N \bigg|\\
&\le \sum_k |\zeta_k^N|  \bigg[ \E Y_{k,2}\left( N^{\gamma}N^{c_k} \int_0^t \lambda_k(\X(s))  -  \lambda_k(\X(s))  \wedge   \lambda_k(\Zm_{\ell}\circ \eta_{\ell}(s))  ds \right)\\
&\hspace{.62in} + \E Y_{k,3}\left(N^{\gamma}N^{c_k}  \int_0^t  \lambda_k(\Zm_{\ell}\circ \eta_{\ell}(s))  - \lambda_k(\X(s))  \wedge  \lambda_k(\Zm_{\ell}\circ \eta_{\ell}(s)) ds  \right)\bigg]\\
&= \sum_k |\zeta_k^N| N^{\gamma}N^{c_k} \int_0^t \E| \lambda_k(X^N(s) )  - \lambda_k(Z^N_{\ell}\circ \eta_{\ell}(s)) | ds\\
&\le N^{\gamma}C  \int_0^t \E |X^N(s) -  Z^N_{\ell}\circ \eta_{\ell}(s) | ds,
\end{align*}
where $C>0$ is some constant and we used that the $\lambda_k$ are assumed to be Lipschitz.
Adding and subtracting the obvious terms yields 
\begin{align}
\begin{split}
 \E | X^N(t) - Z^N_{\ell}(t)| &\le  N^{\gamma} C \int_0^t \E |Z^N_{\ell}(s) - Z^N_{\ell}\circ \eta_{\ell}(s)|ds
    +N^{\gamma} C \int_0^t \E | X^N(s)  - Z^N_{\ell}(s)| ds.
 \end{split}
 \label{eq:bound_ends_exact}
\end{align}
The integrand of the first term on the right hand side of \eqref{eq:bound_ends_exact} satisfies
\begin{align}
     \E |Z^N_{\ell}(s) - Z^N_{\ell}\circ \eta_{\ell}(s)|  &\le \sum_k |\zeta_k^N|N^{\gamma}N^{c_k} \E \int_{\eta_{\ell}(s)}^s \lambda_k(\Zm_{\ell}(\eta_{\ell}(r)) dr    \le \tilde C N^{\gamma} h_{\ell},
     \label{eq:Z_bound}
\end{align}
where $\tilde C>0$ is a constant, and we recall that $\lambda$ is $O(1)$ in our region of interest. Collecting the above yields 
\begin{align*}
	\E | X^N(t) - Z^N_{\ell}(t)  | \le 
 {\widehat C}_1 N^{2\gamma}t h_{\ell} + \widehat{C}_2 N^{\gamma}  \int_0^t \E |X^N(s) - Z^N_{\ell}(s)| ds,
\end{align*}
for some positive constants $\widehat{C}_1,\widehat{C}_2$ that are independent of $N$, $\gamma$, and $T$. The result now follows from Gronwall's inequality.
   \end{proof} 
    
    We note that Lemma \ref{lem:first_bound_exact} is a worst case scenario due to the appearance of the term $N^{\gamma}$ in the exponent.  However, considering the network $S_1\overset{N^{\gamma}}{\to} 2S_1$ (exponential growth), shows this to  be a sharp estimate.  A future research direction will be classifying those networks for which this upper bound can be decreased substantially.  
    
    We are now in position to prove Theorem \ref{thm:MLMC_exact}.

\begin{proof}(of Theorem \ref{thm:MLMC_exact}.)
  We have
\begin{align*}
 X^N(t) - Z^N_{\ell}(t)  =  M^N(t) + \int_0^t F^N(\X(s)) - F^N(\Zm_{\ell}\circ \eta_{\ell}(s)) ds,
\end{align*}
where
\begin{align*}
	M^N(t) &\eqdef  \sum_k \bigg[ Y_{k,2}\left( N^{\gamma}N^{c_k} \int_0^t \lambda_k(\X(s))  -  \lambda_k(\X(s))  \wedge   \lambda_k(\Zm_{\ell}\circ \eta_{\ell}(s))  ds \right)  \\
	&\hspace{.6in}  -   N^{\gamma}N^{c_k} \int_0^t \lambda_k(\X(s))  -  \lambda_k(\X(s))  \wedge   \lambda_k(\Zm_{\ell}\circ \eta_{\ell}(s))  ds  \bigg] \zeta_k^N \\
	&- \sum_k \bigg[ Y_{k,3}\left(N^{\gamma}N^{c_k}  \int_0^t  \lambda_k(\Zm_{\ell}\circ \eta_{\ell}(s))  - \lambda_k(\X(s))  \wedge  \lambda_k(\Zm_{\ell}\circ \eta_{\ell}(s)) ds  \right) \\
	&\hspace{.6in} + N^{\gamma}N^{c_k} \int_0^t \lambda_k(\Zm_{\ell}\circ \eta_{\ell}(s))  - \lambda_k(\X(s))  \wedge  \lambda_k(\Zm_{\ell}\circ \eta_{\ell}(s)) ds \bigg]\zeta_k^N ,
\end{align*}
is a martingale, and
\begin{equation*}
F^N(x) = \sum_k N^{\gamma}N^{c_k} \lambda_k(x) \zeta_k^N.
\end{equation*}
Note that based upon our assumptions, we have that 
\begin{equation}
  |F^N(x) - F^N(y)| \le C N^{\gamma}|x-y|,
  \label{eq:Lipschitz_F}
\end{equation}
 where $C>0$ is a constant that does not depend upon $N$ or $\gamma$.
  The quadratic covariation matrix of $M^N$ is 
\begin{align*}
	[M^N](t) &=   \sum_k \zeta_k^N(\zeta_k^N)^T (J_{k,2}^N(t) + J_{k,3}^N(t)),
\end{align*}
where
\begin{align*}
J_{k,2}^N(t) &\eqdef Y_{k,2}\left( N^{\gamma}N^{c_k} \int_0^t \lambda_k(\X(s))  -  \lambda_k(\X(s))  \wedge   \lambda_k(\Zm_{\ell}\circ \eta_{\ell}(s))  ds \right) \\
J_{k,3}^N(t) &\eqdef Y_{k,3}\left(N^{\gamma}N^{c_k}  \int_0^t  \lambda_k(\Zm_{\ell}\circ \eta_{\ell}(s))  - \lambda_k(\X(s))  \wedge  \lambda_k(\Zm_{\ell}\circ \eta_{\ell}(s)) ds  \right) .
\end{align*}
Thus, 
\begin{equation*}
 \E[M^N](t) =  \sum_k\zeta_k^N(\zeta_k^N)^T N^{\gamma}N^{c_k} \E \int_0^t \big|\lambda_k(X^N(s) )  -  \lambda_k(Z^N_{\ell}\circ \eta_{\ell}(s))  \big| \ ds,
\end{equation*}
and, in particular,
\begin{align}
	\E[M^N]_{ii}(t) = \sum_k (\zeta_{ik}^N)^2N^{\gamma}N^{c_k}   \E \int_0^t \big|\lambda_k(X^N(s) )  -  \lambda_k(Z^N_{\ell}\circ \eta_{\ell}(s))  \big| \ ds.
	\label{eq:qv2_2}
\end{align}
We note that 
\begin{align}
	|X^N(t) -  Z^N_{\ell}(t) |^2  &\le   2   |M^N(t)|^2  + 2   \left| \int_0^t F^N(\X(s)) - F^N(Z^N_{\ell}\circ \eta_{\ell}(s)) ds \right|^2,
	\label{eq:two_terms}
\end{align}
and we may handle the two terms on the right hand side of the above equation separately.

First, by \eqref{eq:qv2_2} and the Burkholder-Davis-Gundy inequality,
\begin{align}
\begin{split}
 \E[ |M^N(t)|^2 ]  &\le \sum_{i} \sum_k (\zeta_{ik}^N)^2N^{\gamma}N^{c_k} \E \int_0^t  \big|\lambda_k(X^N(s) )  -  \lambda_k(Z^N_{\ell}\circ \eta_{\ell}(s))  \big| \ ds \\
&=   \sum_k |\zeta_{k}^N|^2N^{\gamma}N^{c_k} \E \int_0^t  \big|\lambda_k(X^N(s) )  -  \lambda_k(Z^N_{\ell}\circ \eta_{\ell}(s))  \big|  \ ds \\
&\le 2 C  N^{\gamma} N^{-\rho} \E \int_0^t \big| X^N(s) - Z^N_{\ell}\circ \eta_{\ell}(s) \big| \ ds,
\end{split}
\label{eq:last_2}
\end{align}
where $C$ is a constant independent of $N$, $t$, and $\gamma$.
 After adding and subtracting $Z^N_{\ell}(s)$, using \eqref{eq:Z_bound}, and applying Lemma \ref{lem:first_bound_exact}, we conclude that for $t \le T$
 \begin{equation}
 	\E [ |M^N(t)|^2 ]  \le  (c_1 N^{\gamma}T e^{c_2 N^{\gamma} T} ) N^{-\rho} (N^{\gamma} h_{\ell}),
	\label{eq:MN}
 \end{equation}
 for some constants $c_1,c_2$ that do not depend upon $T$, $\gamma$, or $N$, and which will change during the course of the proof.
  
Turning to the second term on the right hand side of \eqref{eq:two_terms}, making use of \eqref{eq:Lipschitz_F} we have for some $C>0$ independent of $T$, $\gamma$, and $N$,
\begin{align}
   \E  \bigg( \int_0^t & |F^N(X^N(s)) - F^N(Z^N_{\ell} \circ \eta_{\ell}(s) ) | ds \bigg)^2\notag \\
     &\le C N^{2\gamma} \E\bigg( \int_0^t  |Z^N_{\ell}\circ \eta_{\ell}(s) - Z^N_{\ell}(s)|ds\bigg)^2   +C N^{2\gamma} \E\bigg( \int_0^t  |X^N(s) -  Z^N_{\ell}(s)|ds\bigg)^2.
   \label{eq:three_terms_2}
\end{align}
The expected value in the first term on the right hand side of \eqref{eq:three_terms_2} can be bounded via
\begin{align}
\begin{split}
	 \E \bigg( \int_0^T  |Z^N_{\ell}\circ \eta_{\ell}(s) - \Zm_{\ell}(s)| ds\bigg)^2 & \le T  \E \int_0^T |Z^N_{\ell}\circ \eta_{\ell}(s) - Z^N_{\ell}(s) |^2 ds\\
	 &= T \sum_{i = 1}^n \int_{t_i}^{t_i + h} \E|Z^N_{\ell}\circ \eta_{\ell}(s) - Z^N_{\ell}(s) |^2 ds.
	\end{split}
	\label{eq:bound1_2}
\end{align}
We have that
\begin{align}
\begin{split}
     \E |Z^N_{\ell}\circ \eta_{\ell}(s) - Z^N_{\ell}(s) |^2  &\le \sum_k |\zeta_k^N|^2\bigg[ N^{\gamma}N^{c_k} \E \int_{\eta_{\ell}(s)}^s \lambda_k(Z^N_{\ell} \circ \eta_{\ell}(r)) dr \\
     &\hspace{.3in}+ N^{2\gamma}N^{2c_k} \E \bigg(\int_{\eta_{\ell}(s)}^s \lambda_k(Z^N_{\ell}\circ \eta_{\ell}(r))  dr\bigg)^2\bigg]\\
     &\le C N^{\gamma} N^{-\rho} h_{\ell} + C N^{2\gamma} h_{\ell}^2,
     \end{split}
     \label{eq:bound1_3}
\end{align}
for some constant $C>0$. Combining \eqref{eq:bound1_2} and \eqref{eq:bound1_3} shows
\begin{equation}
	\E \bigg( \int_0^T  |Z^N_{\ell}\circ \eta_{\ell}(s) - \Zm_{\ell}(s)| ds\bigg)^2 \le  T^2(C N^{\gamma} N^{-\rho} h_{\ell} + C N^{2\gamma} h_{\ell}^2).
	\label{eq:bound_final1}
\end{equation}
 Combining \eqref{eq:bound_final1} with \eqref{eq:three_terms_2} then yields
\begin{align}
\label{eq:final?}
\begin{split}
 \E  \bigg( \int_0^t |F^N(X^N(s)) - F^N(Z^N_{\ell} \circ \eta_{\ell}(s) ) | ds \bigg)^2 \le & \ c_1 N^{2\gamma} T^2 N^{-\rho} (N^{\gamma} h_{\ell})+ c_2 T^2 N^{2\gamma} (N^{\gamma} h_{\ell})^2\\
  &+ c_3t N^{2\gamma}\int_0^t \E |X^N(s) - Z^N_{\ell}(s)|^2 ds,
  \end{split}
\end{align}
for some constants $c_1,c_2,c_3$ that do not depend upon $T$, $N$, or $\gamma$.
Equations \eqref{eq:two_terms}, \eqref{eq:MN}, and \eqref{eq:final?} yield
\begin{align*}
\E[|X^N(t) -  \Zm_{\ell}(t)|^2 ] &\le (c_1 N^{\gamma}T e^{c_2 N^{\gamma} T} ) N^{-\rho} (N^{\gamma} h_{\ell}) + c_1 N^{2 \gamma} T^2  N^{-\rho} (N^{\gamma} h_{\ell})+ c_2 T^2 N^{2\gamma} (N^{\gamma} h_{\ell})^2\\
  &\hspace{.3in}+ c_3t N^{2\gamma}\int_0^t \E |X^N(s) - Z^N_{\ell}(s)|^2 ds.
\end{align*}
The result now follows from Gronwall's inequality.
\end{proof}

 We turn our focus to the proof of Theorem \ref{thm:MLMC}, which is restated here for completeness.\vspace{.1in}
 
\noindent \textbf{Theorem 1.} \textit{
    	Suppose $(\Zm_{\ell},\Zm_{\ell-1})$ satisfy \eqref{eq:Z1} and \eqref{eq:Z2} with $\Zm_{\ell}(0)=\Zm_{\ell-1}(0)$.  Then, there exist functions $C_1,C_2$, that do not depend on $h_{\ell}$, such that 
   \begin{equation*}
      \sup_{t \le T} \E| \Zm_{\ell}(t) - \Zm_{\ell - 1}(t) |^2 \le C_1(N^{\gamma}T)N^{-\rho}(N^{\gamma} h_{\ell}) + C_2(N^{\gamma}T)   (N^{\gamma} h_{\ell})^2.
         \end{equation*} 
In particular, for $\gamma \le 0$ the values 
$
C_1(N^{\gamma}T)
$ 
and
$
C_2(N^{\gamma}T)
$ 
may be bounded above uniformly in $N$.
}

\begin{proof}(of Theorem \ref{thm:MLMC}.)
	A direct proof can be written along the lines of that for Theorem \ref{thm:MLMC_exact}.  A separate, cruder, proof would simply add and subtract $X^N(t)$ to $|Z_{\ell}^N(t) - Z^N_{\ell-1}(t)|^2$ and use Theorem \ref{thm:MLMC_exact} combined with the triangle inequality.
\end{proof}

\section{Implementation issues}
\label{sec:implementation}

The analysis in Sections \ref{sec:new_methods} and \ref{sec:analysis} specified an order of magnitude for the number of paths, $n_{\ell}$, to be used at each level so as to attain the desired accuracy.  This was needed to prove that the computational complexity can be greatly reduced with an appropriate choice of the $n_{\ell}$.  However, the analysis does not tell us what the $n_{\ell}$ should be with precision, nor does it tell us that these are the optimal $n_{\ell}$, which, of course, will depend on the function $f$, and the model itself.  

Letting $V_{\ell}$ denote the variance of $\Ql_{\ell}$ for a given $n_{\ell}$, and $CPU_{{\ell}}$ be the CPU time needed to generate $n_{\ell}$ paths, we know that 
\begin{equation*}
	CPU_{{\ell}} \approx \frac{K_{\ell}}{V_{{\ell}}},
\end{equation*}
for some $K_{\ell}$ as both $CPU_{\ell}$ and $1/V_{\ell}$ scale linearly with $n_{\ell}$.  Further, for a given tolerance, $\epsilon$, we need
\begin{equation}
	\mathsf{Var}(\Ql) = \sum_{\ell} V_{{\ell}} = (\epsilon/1.96)^2,
	\label{eq:constraint}
\end{equation}
for, say, a 95\% confidence interval (where the quantity 1.96 will be changed depending upon the size of the confidence interval desired).  We may approximate each $K_{\ell}$ with a number of preliminary simulations (not used in the full implementation), and then minimize
\begin{equation*}
	\sum_{\ell} \frac{K_{\ell}}{V_{\ell}}, 
\end{equation*}
subject to the constraint \eqref{eq:constraint}.  This will give us target variances, $V_{\ell}$, for each level, and an estimate on the time needed until the computation is completed.  We may then simulate each level until enough paths have been generated for the variance of the estimator at that level to be below the target $V_{\ell}$.  Note that this is similar to the strategy proposed in \cite{Giles2008}.

We make the important observation that with such an optimizing pre-computation, we can choose both $L$ and $\ell_0$, that is the finest and crudest levels, before attempting the full calculation.  This reduces the probability of the approximate processes becoming negative, or non-physical in other ways, during the course of a simulation.  This, in turn, helps keep the approximate processes stable, which leads to increased efficiency.  Further, if we find during such a pre-computation that no choice of levels and number of paths will be significantly faster than using an exact algorithm combined with crude Monte Carlo, then we should simply revert to solely using an exact algorithm.  We may conclude, therefore, that {\em the developed method will never, for any example, be appreciably slower than using an exact algorithm with crude Monte Carlo}.  As will be demonstrated in the next section, however, the method will often times, even in cases not yet predicted by the analysis, be {\em significantly} faster.

Finally, we note that in each of the examples in Section \ref{sec:examples}, 
and each method tested, we use Matlab's built in Poisson random number generator.  Further, we produce the necessary approximate paths in batches ranging from the 100s to 10s of thousands so as to reduce the number of separate calls to the Poisson random number generator.

\section{Examples}
\label{sec:examples}

We present three examples to demonstrate the performance of the proposed method.\vspace{.1in}

\noindent \textbf{Example.}
 We begin by considering a  model of gene transcription and translation also used in \cite{AndMaso2011}:
 \begin{center}
 	$\displaystyle G \overset{25}{\to} G + M, \quad M \overset{1000}{\to} M + P, \quad P+P \overset{0.001}{\to} D, \quad M \overset{0.1}{\to} \emptyset, \quad P \overset{1}{\to} \emptyset. $
 \end{center}
Here, a single gene is being transcribed into mRNA, which is then being translated into proteins, and finally the proteins produce stable dimers.  The final two reactions represent degradation of mRNA and proteins, respectively.  We suppose the system starts with one gene and no other molecules, so $X(0) = (1,0,0)$ where $X_1,X_2,X_3$ give the molecular counts of the mRNA, proteins, and dimers, respectively.  Finally, we suppose that we want to estimate the expected number of dimers at time $T=1$ to an accuracy of $\pm 1$ with 95\% confidence.  Thus, we want the variance of our estimator to be smaller than $(1/1.96)^2 \approx .2603$.  We will also estimate the second moment of the number of dimers, which could be used in conjunction with the mean to estimate the variance.  For comparison purposes, we will use each method discussed in this paper to approximate the mean, and will use an exact method combined with crude Monte Carlo and the unbiased MLMC method to approximate the second moment.
 
 While $\epsilon = 1$ for the unscaled version of this problem, the simulation of just a few paths of the system  shows that there will be approximately 23 mRNA molecules, 3,000 proteins, and 3,500 dimers at time $T=1$.  Therefore, for the scaled system, we are asking for an accuracy of $\widetilde \epsilon = 1/3500 \approx 0.0002857$.  Also, a few paths (100 is sufficient) shows that the order of magnitude of the variance of the normalized number of dimers is approximately 0.11.   Thus, the approximate number of exact sample paths we will need to generate can be found by solving
 \begin{equation*}
 	\frac{1}{n} \mathsf{Var}(\text{normalized \# dimers}) = (\widetilde \epsilon/1.96)^2 \implies n = 5.18 \times 10^6.
 \end{equation*}
 Therefore, we will need approximately five million independent sample paths generated via an exact algorithm.
 
 We also note that with the rough orders of magnitude computed above for the different molecular counts at time $T=1$, we have $N \approx$ 3,500, $\alpha_1 \approx .38$, and  $\alpha_2 = \alpha_3 \approx 1$.  Therefore, we have that $N^\gamma \approx 23,000/3,000 = 7.6 \implies \gamma\approx 0.2485$ for this problem (where we chose the ``stiffest'' reaction for this calculation, which is that of $M \to M + P$).  However, we note that the parameter $\gamma$ changes throughout the simulation and is quite a bit higher near $t \approx 0$.

   Implementing the modified next reaction method, which produces exact sample paths \cite{Anderson2007a}, on our machine\footnote{We used an Apple machine with 8GB Ram and an i7 chip.}  (using Matlab), each path takes approximately 0.03 CPU seconds to generate.  Therefore, the approximate amount of time to solve this particular problem will be 155,000 CPU S, which is about forty three hours.  The outcome of such a simulation is detailed in  Table \ref{table:exact_dimer} where ``\# updates'' refers to the total number, over all paths, of steps, and is used as a crude quantification for the computational complexity of the different methods under consideration.  
 
\begin{table}
\begin{center}
	\begin{tabular}{|c|c|c|c|c|c|}\hline
		 Mean & Variance & $\E X_{3}^2(1)$ & \# paths &  CPU Time &   \# updates    \\  \hline
		 3714.2 $\pm$ 1.0 & $\approx$ 1,232,418 & 1.5035 $\times 10^7$ $ \pm$ 8 $\times 10^3$ & 4,740,000  & 1.49 $\times 10^5$ CPU S   &  8.27 $\times 10^{10}$ \\ \hline
	\end{tabular}
	\end{center}
	\caption{Performance of an exact algorithm with crude Monte Carlo.  The mean number of dimers at time 1 is reported with 95\% a confidence interval.  The approximated variance of the number of dimers is provided for completeness.  An estimate of the second moment is also provided with a 95\% confidence interval.}
	\label{table:exact_dimer}
\end{table}

Next, we solved the problem using Euler tau-leaping with various step-sizes, combined with a crude Monte Carlo estimator.  The results of those simulations are detailed in  Table \ref{table:Euler_dimer}.  Note that the bias of the approximate algorithm has become apparent.\footnote{This data also appears in \cite{AndMaso2011}.}
\begin{table}
\begin{center}
\begin{tabular}{|c|c|c|c|c|c|}\hline
Step-size & Mean  & \# paths &  CPU Time  &  \# updates    \\     \hline
 $h = 3^{-7}$&  3,712.3 $\pm$ 1.0  & 4,750,000 & 13,374.6  S   & $6.2 \times 10^{10}$  \\ [.1 ex] \hline
 $h = 3^{-6}$&  3,707.5 $\pm$ 1.0    &  4,750,000 & 6,207.9 S  &  $2.1 \times 10^{10}$ \\[.1 ex] \hline
 $h = 3^{-5}$&  3,693.4 $\pm$ 1.0   & 4,700,000 & 2,803.9  S &  $6.9 \times 10^9$ \\ [.1 ex]\hline
 $h = 3^{-4}$& 3,655.2 $\pm$ 1.0  & 4,650,000 &  1,219.0  S &  $2.6 \times 10^9$ \\ \hline
\end{tabular}
\end{center}
\caption{Performance of Euler tau-leaping with crude Monte Carlo for the computation of the first moment of $X_3$, the number of dimers. The bias of the method is apparent.  }
\label{table:Euler_dimer}
\end{table}
We then implemented the biased version of MLMC with various step-sizes.  The results of those simulations are detailed in  Table \ref{table:MLMC_dimer_biased}, where the approximations and CPU times should be compared with those of Euler tau-leaping.  The CPU times stated include the time needed to solve the embedded optimization problem discussed in Section \ref{sec:implementation}.  
\begin{table}
\begin{center}
\begin{tabular}{|c|c|c|c|}\hline
Step-size parameters & Mean &   CPU Time &  \# updates    \\     \hline
   $M = 3$, $L = 7$ & 3,712.6 $\pm$ 1.0 & 781.8  S &  1.1 $\times 10^9$  \\ [.1 ex] \hline
   $M = 3$, $L = 6$ & 3,708.5 $\pm$ 1.0  &  623.9  S  &  7.9 $\times 10^8$ \\[.1 ex] \hline
   $M = 3$, $L = 5$ & 3,694.5 $\pm$ 1.0  & 546.9  S &  6.6 $\times 10^8$   \\ [.1 ex]\hline
\end{tabular}
\end{center}
\caption{Performance of biased MLMC with $M = 3$, $\ell_0 = 2$, and $L$ ranging from 7 to 5. The reported times include the time needed for the pre-computations used to choose the number of paths per level as discussed in Section \ref{sec:implementation}.}
\label{table:MLMC_dimer_biased}
\end{table}
Note that the gain in computational complexity, as quantified by the \# updates, over straight tau-leaping with a finest level of $h_L = 3^{-7}$ is 56 fold, with straight tau-leaping taking 17.1 times longer.  Also note that the bias of the approximation method is still apparent.
			  
	Finally, we implemented the unbiased version of MLMC with various step-sizes.  The results of those simulations are detailed in  Table \ref{table:MLMC_unbiased}.  As in the biased MLMC case, the CPU times stated include the time needed to solve the embedded optimization problem discussed in Section \ref{sec:implementation}.
\begin{table}	
\begin{center}
\begin{tabular}{|c|c|r|c|c|}\hline
 Step-size parameters & Mean &  CPU Time & Var. of estimator  & \# updates    \\     \hline
   $M = 3$, $L = 6$ & 3,713.9 $\pm$ 1.0 & 1,063.3  S  & 0.2535 & 1.1 $\times 10^9$  \\ [.1 ex] \hline
  $M = 3$, $L = 5$ & 3,714.7 $\pm$  1.0 & 1,114.9  S  & 0.2565 & 9.4 $\times 10^8$  \\[.1 ex] \hline
   $M = 3$, $L = 4$ & 3,714.2 $\pm$  1.0 & 1,656.6  S  &   0.2595 & 1.0 $\times 10^9$  \\ [.1 ex]\hline
   $M = 4$, $L = 4$ & 3714.2  $\pm$  1.0 &  1,334.8 S & 0.2580  & 1.1  $\times 10^9$  \\ [.1 ex]\hline
   $M = 4$, $L = 5$ & 3,713.8 $\pm$ 1.0 & 1,014.9  S  & 0.2561 & 1.1 $\times 10^9$  \\ \hline
\end{tabular}
\end{center}
\caption{Performance of unbiased MLMC with $\ell_0 = 2$, and $M$ and $L$ detailed above.  The reported times include the time needed for the pre-computations used to choose the number of paths per level as discussed in Section \ref{sec:implementation}.}
\label{table:MLMC_unbiased}
\end{table}
We see that the unbiased MLMC estimator behaves as the analysis predicts: there is no bias for any choice of $M$ or $L$, and the required CPU time are analogous to Euler's method with a course time-step.
Further, the exact algorithm with crude Monte Carlo, by far the most commonly used method in the literature, demanded approximately 80 times more updates and 140 times more CPU time than our unbiased MLMC estimator, with the precise speedups depending upon the choice of $M$ and $L$. 

We feel it is instructive to give more details to at least one choice of $M$ and $L$ for the unbiased MLMC estimator.  For the case with $M=3$, $L=5$, and $\ell_0 = 2$, we provide in Table \ref{table:MoreData35} the relevant data for the different levels.  As already stated in Table \ref{table:MLMC_unbiased},  the total time with the optimization problem was 1,114.9 CPU S, more than the total CPU time reported in Table \ref{table:MoreData35}, which does not include the time needed to solve the optimization problem.  Note that most of the CPU time was taken up at the coarsest level, as is common with MLMC methods and predicted by the analysis.  Also, while the exact algorithm with crude Monte Carlo demanded the generation of almost five million exact sample paths, we needed only 3,900 such paths at our finest level.  This difference is the main reason for the dramatic reduction in CPU time.  Of course, we needed 
more than eight million paths at the coarsest level, but these paths are very cheap to generate.  Finally, we note that the optimization problem divided up the total desired variance into non-uniform sizes with the more computationally intensive levels
generally being allowed to have a higher variance.

In Table \ref{table:MLMC_unbiased_2nd} we provide the data pertaining to the estimate of the second moment using the unbiased version of MLMC with $M = 3$ and $L = 5$
that appears in the second row of Table~\ref{table:MLMC_unbiased}. The 95\% confidence interval for the second moment is 1.5031 $\times 10^7$ $\pm $ 9 $ \times 10^3$, which should be compared with the confidence interval generated using an exact method.

	\begin{table}
\begin{center}
	\begin{tabular}{|c|c|c|c|c|c|}\hline
		Level & \# paths &  Mean &  Var.  estimator & CPU Time  & \# updates    \\\hline
		 $(X,Z_{3^{-5}})$ & 3,900 & 20.1 & 0.0658  & 279.6  S &  6.8 $\times 10^7$  \\ \hline
		  $(Z_{3^{-5}},Z_{3^{-4}})$ & 30,000 & 39.2 &0.0217 & 49.0 S  & 8.8 $\times 10^7$ \\ \hline
		  $(Z_{3^{-4}},Z_{3^{-3}})$ & 150,000 & 117.6 & 0.0179 & 71.7  S    & 1.5  $\times 10^8$  \\ \hline
		  $(Z_{3^{-3}},Z_{3^{-2}})$ & 510,000 & 350.4 & 0.0319  & 112.3  S &  1.7 $\times 10^8$  \\ \hline
		  Euler, $h = 3^{-2}$ & 8,630,000 & 3,187.4 &  0.1192 & 518.4 S   & 4.7  $\times 10^8$ \\ \hline
		  Totals & N.A. &3,714.7 & 0.2565 & 1031.0  S   & 9.5 $\times 10^8$\\ \hline
	\end{tabular}
	\end{center}
	\caption{Details of the different levels for the implementation of the unbiased MLMC method with $M = 3$, $L = 5$, and $\ell_0 = 2$.  By $(X,Z_{3^{-5}})$ we mean the level in which the exact process is coupled to the approximate process with $h = 3^{-5}$, and by $(Z_{3^{-\ell}},Z_{3^{-\ell+1}})$ we mean the level with $Z_{3^{-\ell}}$ coupled to $Z_{3^{-\ell+1}}$.  }
	\label{table:MoreData35}
	\end{table}

	\begin{table}	
\begin{center}
	\begin{tabular}{|c|c|c|c|}\hline
		Level &  Estimate &  Var. of estimator  & 95\% Confidence inteval    \\\hline
		 $(X,Z_{3^{-5}})$ & 157,000  &  5.8 $\times 10^6$  & N.A.   \\ \hline
		  $(Z_{3^{-5}},Z_{3^{-4}})$ & 303,000 & 2.0 $ \times 10^6$ & N.A.  \\ \hline
		  $(Z_{3^{-4}},Z_{3^{-3}})$ & 894,000 & 2.2 $\times 10^6$ & N.A.  \\ \hline
		  $(Z_{3^{-3}},Z_{3^{-2}})$ & 2.4888 $\times 10^6$ & 4.2 $\times 10^6$ &N.A.  \\ \hline
		  Euler, $h = 3^{-2}$ & 1.1188 $\times 10^7$ & 5.7 $\times 10^6$ & N.A. \\ \hline
		  Totals & 1.5031 $\times 10^7$ &  2.0 $\times 10^7$ & 1.5031 $\times 10^7$ $\pm $ 9 $ \times 10^3$ \\ \hline
	\end{tabular}
	\end{center}
	\caption{Details of the different levels for the implementation of the unbiased MLMC method with $M = 3$, $L = 5$, and $\ell_0 = 2$ for the approximation of the second moment.   }
\label{table:MLMC_unbiased_2nd}
\end{table}

\vspace{.2in}

\noindent \textbf{Example.}  We turn now to a simple example 
that allows us to study how the behavior of the developed methods depends upon the parameter $\gamma$.  Consider the family of models indexed by $\theta$,
\begin{align*}
	&A \overset{\theta}{\underset{\theta}{\rightleftarrows}} B,
	\end{align*}
with,
\begin{align*}
	X_A(0) = X_B(0) = \lfloor \text{1,000}\, \theta^{-1}\rfloor,
\end{align*}
where $\lfloor x\rfloor$ is the greatest integer less than or equal to $x$.  The stochastic equation governing $X_A$, giving the number of $A$ molecules, is
\begin{equation*}
	X_A(t) = X_A(0) + Y_1\left( \int_0^t \theta (2,000 \theta^{-1} - X_A(s)) ds\right) - Y_2\left( \int_0^t \theta X_A(s) ds\right),
\end{equation*}
with the Euler approximation given by
\begin{equation*}
	Z_A(t) = Z_A(0) + Y_1\left( \int_0^t \theta (2,000\theta^{-1} - Z_A\circ \eta(s)) ds\right) - Y_2\left( \int_0^t \theta Z_A\circ \eta(s) ds\right).
\end{equation*}

Letting $N = X_A(0)$, we see that $\theta = N^{\gamma}$, implying 
\begin{equation*}
	\gamma = \ln(\theta)/\ln(N).
\end{equation*}
Note, in particular, that  $\gamma > 0$ $\iff$ $\theta > 1$, with $\gamma$ being a strictly increasing function of $\theta$.  Therefore, we may test the dependence of the behavior of the MLMC method on this model by varying the single parameter $\theta$.   We will let $\theta$ range from $0.1$ to $1,000$ and for each $\theta$ use both an exact method with crude Monte Carlo and the unbiased MLMC method developed here to estimate the mean number of $A$ molecules at time 1.  We choose different values for our tolerance parameter, $\epsilon$, for different values of $\theta$.  The results of these computations can be found in Table \ref{table_simp_ex_theta}. 

\begin{figure}
\begin{center}
\includegraphics[height=3in]{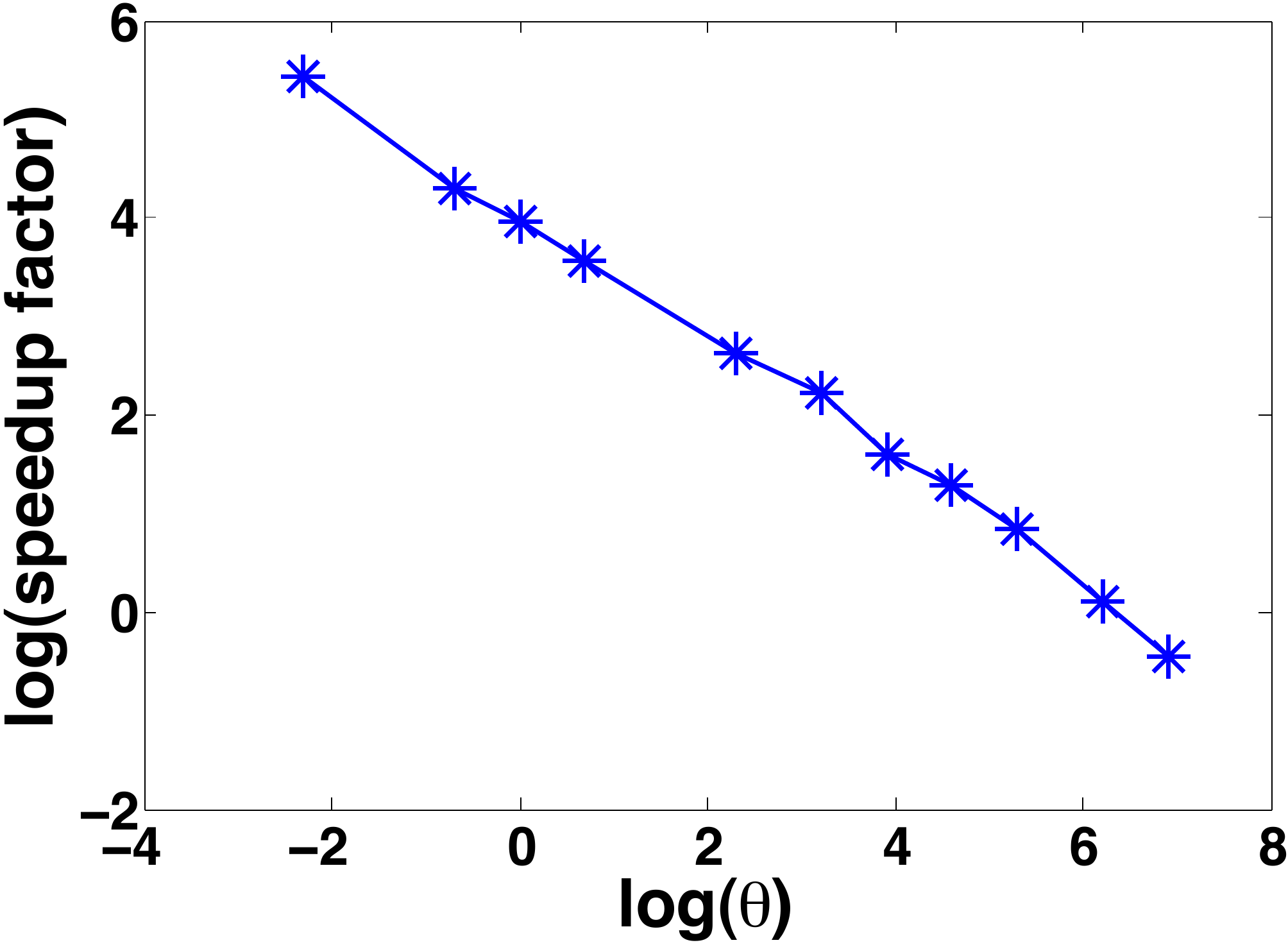}
\end{center}
\caption{Log-log plot for the speedup factor of unbiased MLMC over an exact method.  The best fit line, not shown, is $4 - 0.62x$.}
\label{log-log_speedup}
\end{figure}

  The analysis of this paper predicts that as $\theta$ increases, MLMC should progressively lose its computational advantage over an exact algorithm, and this is borne out in the data provided in Table \ref{table_simp_ex_theta}.  Note, however, that the unbiased version of MLMC remains significantly more efficient than an exact algorithm until $\theta = 1,000$, in which case $X_A(0) = 1$.  Having $\theta = $ 1,000  is arguably the {\em worst} case scenario for an approximate algorithm such as tau-leaping, and the coupling method performs slightly worse than using an exact method with crude Monte Carlo.  As discussed in Section \ref{sec:implementation}, we would normally in this case simply use the exact method alone. However, we report the MLMC data for the sake of comparison.  Further, it is extremely encouraging that there were still large gains in efficiency even when $\theta \in \{25, 50, 100\}$, something not predicted by the current analysis.  Interestingly, the speedup factor of MLMC over an exact method appears, for this example, to be a function of $\theta$.  Specifically, as demonstrated by the log-log plot in Figure \ref{log-log_speedup}, we observed the relation
\begin{equation*}
  	\text{Speedup factor} \approx 54.6\,  \theta^{-0.62}.
  \end{equation*}
\begin{table}
\begin{center}
	\begin{tabular}{|c|c|c|c|c|c|}\hline
		$\theta$ & Method & Estimate of $X_A(1)$ & CPU time & \# Paths & Speedup \\\hline
		.1 & Crude MC & 9,999.97 $\pm$ 0.20 & 1,476.9 S & 164,800 & N.A. \\ \hline
		.1 & MLMC: $M = 3$, $L = 4$, $\ell_0=0$ &  10,000.07 $\pm$  0.19  & 6.5 S & N.A. & 227.2 \\ \hline \hline
		 .5 & Crude MC & 1,999.90 $\pm$ 0.16 & 1,110.9 S & 124,100 & N.A. \\ \hline
		 .5 & MLMC: $M = 3$, $L = 4$, $\ell_0=1$ & 2,000.10  $\pm$  0.16  & 15 S & N.A. & 74.1  \\ \hline \hline
		 1 & Crude MC & 999.96 $\pm$  0.11 & 1,464.4 S & 163,800 & N.A. \\ \hline
		  1 & MLMC: $M = 3$, $L = 5$, $\ell_0=2$ & 999.99 $\pm$  0.11 & 28 S & N.A. & 52.3  \\ \hline \hline
		 2 & Crude MC & 500.01 $\pm$ 0.11 & 739.9 S & 82,700 & N.A.  \\ \hline
		2 & MLMC: $M = 6$, $L = 6$, $\ell_0=3$ & 499.96 $\pm$ 0.11 & 21 & N.A. & 35.2 \\ \hline  \hline
		10 & Crude MC & 99.983 $\pm$ 0.044 & 900.2 S & 100,600 & N.A.  \\ \hline
		10 & MLMC:  $M = 3$, $L = 6$, $\ell_0 = 5$ & 99.965 $\pm$ 0.044 & 65 S & N.A. & 13.8 \\ \hline  \hline
		25 & Crude MC &  40.012 $\pm$ 0.028 & 898.0 S & 100,500 & N.A. \\ \hline
		25 & MLMC: $M = 3$, $L = 6$, $\ell_0 = 6$ & 39.996 $\pm$ 0.028 & 98 S & N.A. & 9.2  \\ \hline \hline
		50 & Crude MC & 20.008 $\pm$  0.0139& 1,789.0 S & 200,200 & N.A. \\ \hline 
		50 & MLMC:  $M = 3$, $L = 7,$ $\ell_0 = 7$ & 20.005 $\pm$ 0.0138  & 360 S & N.A. & 5.0 \\ \hline \hline
		100 & Crude MC & 10.002 $\pm$ 0.0139  & 892.6 S & 100,100 & N.A. \\ \hline
		100 & MLMC:  $M = 3$, $L = 7$, $\ell_0 = 7$ &  9.988 $\pm$ 0.0138 & 250 S & N.A. &  3.6  \\ \hline \hline
		200 & Crude MC & 4.9996 $\pm$  0.0088 & 1,120.3 S &125,400 & N.A. \\ \hline
		200 & MLMC: $M = 3$, $L = 7$, $\ell_0 = 7$ & 4.9958 $\pm$ 0.0087  &  486 S & N.A. &2.3 \\ \hline \hline
		500 & Crude MC & 2.0029 $\pm$  0.0044 & 1,781.6 S & 199,400 & N.A. \\ \hline
		500 & MLMC: $M = 3$, $L = 7$, $\ell_0 = 7$  & 1.9953 $\pm$ 0.0044 & 1,625.9 S & N.A. &    1.1 \\ \hline \hline
		1,000 & Crude MC & 1.0038 $\pm$ 0.0043 & 897.2 & 100,200& N.A. \\ \hline
		1,000 & MLMC: $M = 3$, $L= 7$, $\ell_0 = 7$ & 1.0015 $\pm$  0.0044 & 1,412.3 S &N.A. &  0.64 \\ \hline
	\end{tabular}
\end{center}
\caption{Approximation of $X_A(1)$ with 95\% confidence intervals.  Note that the speedup factor decreases as $\theta$ increases, with MLMC becoming less efficient than an exact method when $\theta = 1,000$.}
\label{table_simp_ex_theta}
\end{table}

\vspace{.2in}
 
  \noindent \textbf{Example.}  We finish with a model of viral kinetics first developed in \cite{Yin2002} by  Yin et al., and subsequently studied by Haseltine and Rawlings in \cite{Haseltine2002}, Ball et al., in \cite{Ball06}, and E et al., in \cite{EVE2007}.  One reason the  interest in this model has been so high from the biological, engineering, and  mathematical communities is that it exemplifies a feature of many stochastic models arising in the biosciences: a separation of time scales.  We will use this model to demonstrate that one of the main ideas of this paper, the coupling, is not restricted to the use of approximate methods defined using time-discretizations.
  
   The model includes four time-varying ``species'': the viral genome ($G$), the viral structural protein ($S$), the viral template ($T$), and the secreted virus itself ($V$). We denote these as species 1, 2, 3, and 4, respectively, and let $X_i(t)$ denote the number of molecules of species $i$ in the system at time $t$. The model involves six reactions, 
\begin{align*}
   	R_1:&\quad  T \overset{\kappa_1}{\to} T + G, \hspace{1in} \kappa_1 = 1, \\[1ex]
	R_2:& \quad G \overset{\kappa_2}{\to} T,  \hspace{1.3in} \kappa_2 = 0.025,  \\[1ex]
	R_3:& \quad T \overset{\kappa_3}{\to} T+S, \hspace{1.03in} \kappa_3 = 1000, \\[1ex]
	R_4:&\quad  T \overset{\kappa_4}{\to} \emptyset, \hspace{1.35in} \kappa_4 = 0.25, \\[1ex]
	R_5:& \quad S \overset{\kappa_5}{\to} \emptyset, \hspace{1.37in} \kappa_5 = 2, \\[1ex]
	R_6:&\quad G + S \overset{\kappa_6}{\to} V, \hspace{1.06in} \kappa_6 = 7.5 \times 10^{-6},
\end{align*}
where the units of time are in days.
The stochastic equations  for this model are
\begin{align*}
	X_1(t) &= X_1(0) + Y_1\left( \int_0^t X_3(s)ds\right) - Y_2\left( 0.025 \int_0^t X_1(s)ds\right)\\
	&\hspace{.4in} - Y_6\left( 7.5 \times 10^{-6} \int_0^t X_1(s) X_2(s) ds\right)\\
	X_2(t) &= X_2(0) + Y_3\left( 1000 \int_0^t X_3(s)ds\right) - Y_5\left( 2 \int_0^t X_2(s)ds\right)\\
	&\hspace{.4in} - Y_6\left( 7.5 \times 10^{-6} \int_0^t X_1(s) X_2(s) ds\right)\\
	X_3(t) &= X_3(0) + Y_2\left( 0.025 \int_0^t X_1(s)ds\right)- Y_4\left( 0.25 \int_0^t X_3(s)ds\right)\\
	X_4(t) &= X_4(0) + Y_6\left( 7.5 \times 10^{-6} \int_0^t X_1(s) X_2(s) ds\right).
\end{align*}	

Following \cite{EVE2007}, we assume an initial condition of $X(0) = (0,0,10,0)^t$.  We see that whenever  the number of viral templates is positive, that is whenever $X_3 >0$, the rates of the third and fifth reactions will be substantially larger than the others.  At the times when $X_3>0$ and $X_2 = O(1)$, we have that $\gamma \gg 1$, with $\gamma$ remaining large until $X_2 = O(1000)$.  However, even when $X_2 = O(1000)$, the natural time-scale of the problem is $O(1/1000)$, whereas the time-scale in which we would like to answer questions  is $O(1)$.  

Instead of implementing our MLMC method directly, we take an alternative approach that makes use of the idea of the coupling, though not the multi-level aspect of the paper.   That is, we will build an approximate process $Z$ that will be used as a control variate for $X$.  Towards that end, note that when the number of templates is positive, reactions 3 and 5 are much faster than the others.  Ignoring the other reactions, we see that when $X_3 >0$, the ``system'' governing the dynamical behavior of $S$ is 
 \begin{align*}
 	\emptyset \overset{1000X_3(t)}{\underset{2}{\rightleftarrows}}& S,
 \end{align*}
 which has an equilibrium distribution that is Poisson with a parameter of $500X_3(t)$, see \cite{AndCraciunKurtz}.  
  Believing that we may use this mean value of $X_2(s)$ in the integrated intensity of reaction 6, that is
 \begin{equation}
 	\int_0^t X_1(s) X_2(s) ds \approx \int_0^t X_1(s) (500X_3(s)) ds,
	\label{viral:approx}
 \end{equation}
 we hope a good approximate model for $G$, $T$, and $V$, which we denote by $Z = (Z_1,Z_3,Z_4)$ so as to remain consistent with the enumeration of $X$, is
 \begin{align}
 \begin{split}
	Z_1(t) &= X_1(0) + Y_1\left( \int_0^t Z_3(s)ds\right) - Y_2\left( 0.025 \int_0^t Z_1(s)ds\right)\\
	&\hspace{.4in} - Y_6\left( 3.75 \times 10^{-3} \int_0^t Z_1(s) Z_3(s) ds\right)\\
	Z_3(t) &= X_3(0) + Y_2\left( 0.025 \int_0^t Z_1(s)ds\right)- Y_4\left( 0.25 \int_0^t Z_3(s)ds\right)\\
	Z_4(t) &= X_4(0) + Y_6\left( 3.75 \times 10^{-3} \int_0^t Z_1(s) Z_3(s) ds\right).
	\end{split}
	\label{eq:viral_Z}
\end{align}
Note that while $Z$ is an approximate model of $X$, it is still a valid continuous time Markov chain satisfying the natural non-negativity constraints.  In particular, there is no time-discretization parameter in $Z$, which is where many technical problems related to tau-leaping (stability concerns, negativity of molecular counts, etc.) arise. 

We will now couple the two processes in a manner similar to \eqref{eq:unscaled_coupled_exact} and build our estimator.  Let $\zeta_k$ denote the reaction vector for the $k$th reaction.  Let $\lambda_6(X) = 7.5 \times 10^{-6} X_1X_2$, and $\Lambda_6(Z) = 3.75 \times 10^{-3}Z_1Z_3$.  For arbitrary $f$, we can  estimate $\E f(X(T))$ via
\begin{align}\label{ex:control}
	\E f(X(T)) = \E (f(X(T)) - f(Z(T))) + \E f(Z(T)),
\end{align}
where $\E f(Z(T))$ is estimated by crude Monte Carlo using the representation \eqref{eq:viral_Z},
which is relatively cheap to simulate, and we estimate $\E (f(X(T)) - f(Z(T)))$ using independent realizations from the coupled processes $(X,Z)$ below
{\small
\begin{align}
\begin{split}
	X(t) &= X(0) + Y_{1,1}\left(\int_0^t \min\{X_3(s),Z_3(s)\} ds\right)\zeta_1 + Y_{1,2}\left(\int_0^t X_3(s) - \min\{X_3(s),Z_3(s)\} ds\right)\zeta_1\\
	& + Y_{2,1}\left(0.025 \int_0^t \min\{X_1(s),Z_1(s)\} ds\right)\zeta_2 + Y_{2,2}\left(0.025\int_0^t X_1(s) - \min\{X_1(s),Z_1(s)\} ds\right)\zeta_2\\
	& + Y_3\left( 1000 \int_0^t X_3(s)ds\right)\zeta_3\\
	& +Y_{4,1}\left( 0.25 \int_0^t \min\{X_3(s),Z_3(s)\}(s)ds\right)\zeta_4 + Y_{4,2}\left( 0.25 \int_0^tX_3(s) -  \min\{X_3(s),Z_3(s)\}(s)ds\right)\zeta_4\\
	& + Y_5\left( 2 \int_0^t X_2(s)ds\right)\zeta_5\\
	& + Y_{6,1}\left(  \int_0^t \min\{\lambda_6(X(s)),\Lambda_6(Z(s))\} ds\right)\zeta_6 - Y_{6,2}\left(  \int_0^t \lambda_6(X(s)) - \min\{\lambda_6(X(s)), \Lambda_6(Z(s))\} ds\right)\zeta_6\\
	Z(t) &=  Y_{1,1}\left(\int_0^t \min\{X_3(s),Z_3(s)\} ds\right)\zeta_1 + Y_{1,3}\left(\int_0^t Z_3(s) - \min\{X_3(s),Z_3(s)\} ds\right)\zeta_1\\
	& + Y_{2,1}\left(0.025 \int_0^t \min\{X_1(s),Z_1(s)\} ds\right)\zeta_2 + Y_{2,3}\left(0.025\int_0^t Z_1(s) - \min\{X_1(s),Z_1(s)\} ds\right)\zeta_2\\
	& +Y_{4,1}\left( 0.25 \int_0^t \min\{X_3(s),Z_3(s)\}(s)ds\right)\zeta_4 + Y_{4,3}\left( 0.25 \int_0^tZ_3(s) -  \min\{X_3(s),Z_3(s)\}(s)ds\right)\zeta_4\\
	& + Y_{6,1}\left(  \int_0^t \min\{\lambda_6(X(s)),\Lambda_6(Z(s))\} ds\right)\zeta_6 - Y_{6,3}\left(  \int_0^t \Lambda_6(Z(s)) - \min\{\lambda_6(X(s)), \Lambda_6(Z(s))\} ds\right)\zeta_6,
\end{split}
\label{major_coupling}
\end{align}}
where the $Y_{k,i}$'s are independent, unit-rate Poisson processes.  Note that we have coupled the process through the reaction channels 1, 2, 4, and 6, in the usual way, though not through 3 or 5, which are not incorporated in the model for $Z$.  Simulation of the coupled processes, which is itself just a continuous time Markov chain in $Z^6_{\ge 0}$, may proceed by any exact algorithm. Here we used the next reaction method \cite{Anderson2007a,Gibson2000}.

Supposing we want to estimate  $\E X_{4}(20)$, giving the mean number of virus molecules at time 20,  we calculate this value using both a naive application of the next reaction method with crude Monte Carlo, and the control variate approach of \eqref{ex:control} with the coupling \eqref{major_coupling}.  The details of the two computations are found in Table \ref{table:viral}.  We see that the crude Monte Carlo implementation required 60 times more updates and 22
times more CPU seconds than the control variate/coupling approach, again demonstrating the usefulness of the core ideas of this paper.

\begin{table}	
\begin{center}		
		\begin{tabular}{|c|c|c|c|c|}\hline
		Method &  Approximation & \# paths &  CPU Time  & \# updates    \\  \hline
		  Crude Monte Carlo & 13.85 $\pm$ 0.07 & 75,000  &  24,800 CPU S & $1.45 \times 10^{10}$  \\\hline
		  Coupling & 13.91 $\pm$  0.07  & N.A. &   1,118.5 CPU S  &  $2.41 \times 10^8$\\ \hline
			\end{tabular}
	\end{center}
	\caption{Details of the approximated expected virus level at time 20 using crude Monte Carlo with an exact algorithm and a control variate approach using the coupling \eqref{major_coupling}.}
\label{table:viral}
\end{table}

\section{Conclusions}
\label{sec:conclusions}

This work focused on the Monte Carlo approach to estimating expected values
of continuous time Markov chains.
In this context there is a trade off between the accuracy and the cost of each Monte Carlo sample.   
Exact samples are available, but these are typically very expensive, especially in our target application of biochemical
kinetics. Approximate samples can be computed by tau-leaping, with the bias governed by 
 a discretization parameter, $h$. 
A realistic analysis of the cost of tau-leaping must acknowledge the importance of 
system scaling.  In particular, for a fixed system, 
in the limit $h \to 0$ tau-leaping becomes infinitely more expensive than
exact sampling, since it needlessly refines the waiting times between reactions.  
In this work, we studied tau-leaping in a general setting that incorporates system scaling without taking asymptotic limits.
Motivated by the work of Giles
\cite{Giles2008} on diffusion processes, we then introduced a new multi-level version of the algorithm that combines coordinated
pairs of tau-leaping paths at different $h$ resolutions.  The two main conceptual advances in this work
were (a) pointing out the practical benefits of a coupling process that had previously been 
introduced solely as a theoretical tool, and (b) exploiting the availability of an exact sampling 
algorithm to give an unbiased estimator.  Our theoretical analysis of the 
computational complexity showed that the new algorithm  
dramatically outperforms the existing state of the art in a wide range of scaling regimes, including 
the classical scaling arising in chemical kinetics. 
The new algorithm is straightforward to summarize and implement, and numerical results 
confirmed that the predicted benefits can be seen in practice.

There are several avenues for future work in this area, including,
\begin{itemize}
 \item using Quasi-Monte Carlo sampling to 
improve practical performance,
 \item customizing the method in the context of multi-scale or hybrid models, where it is possible 
         to exploit special structure in the form of fast/slow reactions or species, or
       where  the discrete space Markov chain is coupled to diffusion or ODE models,    
 \item extending the theoretical analysis to the $\gamma > 0$ regime, in order to explain 
why we continued to observe excellent results in practice,
\item using the coupling idea without discretization to obtain a control variate method
  that exploits specific problem structure, as illustrated in the third example of 
section~\ref{sec:examples}.
\end{itemize} 

\noindent
\textbf{Acknowledgement}\\
The authors are grateful to 
the Banff International Research Station 
for Mathematical Innovation and Discovery for supporting their 
attendance at 
the 
workshop on 
\emph{Multi-scale Stochastic Modeling of Cell Dynamics}, January 2010,
where this collaboration began.
We also thank Mike Giles for very useful feedback on an earlier version of this manuscript.

\bibliographystyle{amsplain} \bibliography{MLMC_TAU}

\end{document}